\documentclass[11pt]{article}   
\usepackage{amssymb,amscd,latexsym} 
\usepackage{graphicx,color}  
\usepackage{amsmath}
\usepackage{amsthm}
\usepackage{mathdots}
\usepackage[all]{xy}
\newcommand{\rar}{\rightarrow}

\newcommand{\llar}{-\kern-5pt-\kern-5pt\longrightarrow}
\newcommand{\surjects}{\twoheadrightarrow}

\newtheorem{Theorem}{Theorem}[section]
\newtheorem{Lemma}[Theorem]{Lemma}

\newtheorem{Proposition}[Theorem]{Proposition}
\newtheorem{Remark}[Theorem]{Remark}

\newtheorem{Conjecture}[Theorem]{Conjecture}

\newtheorem{Question}[Theorem]{Question}
\def\sqr#1#2{{\vcenter{\hrule height.#2pt
        \hbox{\vrule width.#2pt height#1pt \kern#1pt
            \vrule width.#2pt}
        \hrule height.#2pt}}}
\def\phi{\varphi}
\def\demo{\noindent{\bf Proof. }}
\def\square{\mathchoice\sqr64\sqr64\sqr{4}3\sqr{3}3}
\def\qed{\hspace*{\fill} $\square$}

\def\xx{{\bf x}}
\def\yy{{\bf y}}

\def\XX{{\bf X}}

\def\vv{{\bf v}}

\def\codim{{\rm codim}\,}

\def\rk{{\rm rank}\,}

\def\restr{{\kern-1pt\restriction\kern-1pt}}

\def\pp{{\mathbb P}}

\begin{document}
\begin{center}
{\Large{\bf\sc Degenerations of the generic square matrix. \\Polar map and determinantal structure}}
\footnotetext{AMS Mathematics
Subject Classification (2010   Revision). Primary 13C40, 13D02, 14E05, 14M12; Secondary  13H10,14M05.}

\vspace{0.3in}

{\large\sc Rainelly Cunha}\footnote{Under a CAPES Doctoral scholarship. The paper contains parts of this author's ongoing PhD thesis} \quad\quad
{\large\sc Zaqueu Ramos}\footnote{Partially supported by a CNPq post-doctoral fellowship (151229/2014-7)} \quad\quad
 {\large\sc Aron  Simis}\footnote{Partially
 	supported by a CNPq grant (302298/2014-2) and by a CAPES-PVNS Fellowship (5742201241/2016).}

\end{center}

\tableofcontents

\begin{abstract}
One studies certain degenerations of the generic square matrix over a field $k$ along with its main related structures, such as the determinant of the matrix, the ideal generated by its partial derivatives, the polar map defined by these derivatives, the Hessian matrix and the ideal of the submaximal minors of the matrix. The main tool comes from commutative algebra,  with emphasis on ideal theory and syzygy theory. The structure of the polar map is completely identified and the main properties of the ideal of submaximal minors are determined.
Cases where the degenerated determinant has non-vanishing Hessian determinant show that the former is a factor of the latter with the (Segre) expected multiplicity, a result treated by Landsberg-Manivel-Ressayre by geometric means.
Another byproduct is an affirmative answer to a question of F. Russo concerning the codimension in the polar image of the dual variety to a hypersurface.

\end{abstract}

\section*{Introduction}

As implicit in the title, we aim at the study of certain degenerations of the generic square matrix over a field $k$ along with its main related structures.
The degenerations one has in mind will be carried by quite simple homomorphisms of the ground polynomial ring generated by the entries of the matrix over $k$,
mapping any entry to another entry or to zero.
Since the resulting degenerated matrices often ``forget'' their generic origins, the study of the related structures becomes a hard step.

Concretely, let $\mathcal M$ denote an $m\times m$ matrix which is a degeneration of the  $m\times m$ generic matrix and let $R$ denote the polynomial ring over $k$ generated by its entries.
The related structures will mean primevally the determinant $f\in R$ of $\mathcal M$, the corresponding Jacobian ideal $J\subset R$, the Hessian matrix $H(f)$ of $f$, the polar map of $f$ defined by its partial derivatives, and the ideal $I\subset R$ of submaximal minors. 
The approach throughout takes mainly the commutative algebra side of the structures, hence ideal theory and homological aspects play a dominant role. 

The usual geometric way is to look at $J$ as defining the base scheme of the polar map without further details about its scheme nature.
Here, a major point is to first understand the ideal theoretic features of $J$ such as its codimension and some of its associated prime ideals, as well as the impact of its linear syzygies. In a second step one focus on the intertwining between $J$ and the ideal $I\subset R$ of the submaximal minors of the matrix. As per default, we will have $J\subset I$ and typically of the same codimension. Of great interest is to find out when $I$ is a prime ideal and to study its potential nature as the radical of $J$ or of its unmixed part. As far as we know, for the question of primality  there are two basic techniques. One is based on verifying beforehand the normality of $R/I$, the other through a lucky application of the results on $s$-generic matrices, as developed in \cite{Eisenbud1} and \cite{Eisenbud2}. We draw on both, depending on the available players.
A beautiful question, still open in general to our knowledge, is to determine when the singular locus of the determinantal variety defined by $I$ is set-theoretic defined by the immediately lower minors, just as happens in  the generic case.

Another major focus is the structure of the polar map and its image (freely called {\em polar image}) in terms of its birational potentiality. Since the polar map is a rational map of projective space to itself, the image can be defined in terms of the original coordinates (variables of $R$) and we may at times make this abuse.
Its homogeneous coordinate ring is known as the special fiber (or special fiber cone) of the ideal $J$ and plays a central role in the theory of reductions of ideals. Its Krull dimension is also known as the {\em analytic spread} of the ideal $J$. In this regard, we face two basic questions: first, to compute the analytic spread of $J$ (equivalently, in characteristic zero,  the rank of the Hessian matrix $H(f)$ of $f$); second, to decide when $J$ is a minimal reduction of the ideal $I$ of submaximal minors.

Complete answers to the above questions are largely dependent upon the sort of degenerations of the generic square matrix one is looking at.
Throughout this work all matrices will have as entries either variables in a polynomial ring over a field or zeros, viewed as degenerations of the generic square matrix.
The prevailing tone of this study is to understand the effect of such degenerations on the properties of the underlying ideal theoretic structures.
By and large, the typical degeneration one has in mind consists in replacing some of the entries of the generic matrix by way of applying a homomorphism of the ambient polynomial ring to itself.
There will be noted differences regarding the behavior of a few properties and numerical invariants, such as codimension, primality, Gorensteiness and Cohen--Macaulayness. 
In particular, degenerating variables to zero is a delicate matter and may depend on strategically located zeros. 

In this work we focus on two basic instances of these situations. 

The first, called informally ``entry cloning'' is dealt with in Section~\ref{CLONING}. Here we show that $f$ is homaloidal. For this, we first prove that the Jacobian ideal $J$ has maximal linear rank and that the determinant of $H(f)$ does not vanish, both requiring quite some {\em tour de force}.
By using the birationality criterion in Theorem~\ref{basic_criterion}, it follows that the polar map is a Cremona map -- i.e., $f$ is a homaloidal polynomial.

We then move on to the ideal $I$ of the submaximal minors. It will be a Gorenstein ideal of codimension $4$, a fairly immediate consequence of specialization. Showing in addition that it is a prime ideal required a result of Eisenbud drawn upon the $2$-generic property of the generic matrix -- we believe that $R/I$ is actually a normal ring.
It turns out that $I$ is the minimal primary component of $J$ and the latter defines a double structure on the variety $V(I)$ with a unique embedded component, the latter being a linear subspace of codimension $4m-5$.
An additional result is that the rational map defined by the submaximal minors is birational onto its image. We give the explicit form of the image, through its defining equation, a determinantal expression of degree $m-1$.
From the purely algebraic side, this reflects on showing that the ideal $J$ is not a reduction of its minimal component $I$.

The last topic of the section is the structure of the dual variety $V(f)^*$ of $V(f)$. Here we show that $V(f)^*$ is an arithmetically Cohen--Macaulay variety that has very nearly the structure of a ladder determinantal ideal of $2$-minors and has dimension $2m-2$.
We close the section by proving that $f$ divides its Hessian determinant and has the {\em expected multiplicity} $m(m-2)-1$ thereof in the sense of B. Segre. 

The second alternative is dealt with in Section~\ref{zeros}.
We replace generic entries by zeros in a strategic position to be explained in the text. For any given $1\leq r\leq m-2$, the degenerated matrix will acquire ${r+1 \choose 2}$ zeros.
We prove that the ideal $J$ still has maximal linear rank. This time around, the Hessian determinant vanishes and the image of the polar map is shown to have dimension $m^2-r(r+1)-1$. Moreover, its homogeneous coordinate ring is a ladder determinantal Gorenstein ring.

Moving over again to the dual variety $V(f)^*$ we find that it is a ladder determinantal variety of dimension $2m-2$ defined by $2$-minors. Thus, it has codimension $(m-1)^2-r(r+1)$ in the polar image. This result answers a question (oral communication) of F. Russo as to whether there are natural examples where the codimension of the dual variety in the polar image is larger than $1$, assuming that the Hessian determinant of $f$ vanishes -- here the gap is actually arbitrarily large and in addition it is representative of a well structured class of determinantal hypersurfaces.
We note that $V(f)^*$ is in particular an arithmetically Cohen--Macaulay variety.
It is arithmetically Gorenstein if  $r=m-2$.

In the sequel, as in the previous section, our drive is the nature of the ideal $I$ of submaximal minors. Once again, we have geometry and algebra. The main geometric result is that these minors define a birational map onto its image and the latter is a cone over the polar variety of $f$ with vertex cut by ${{r+1}\choose {2}}$ coordinate hyperplanes.
The algebraic results are deeper in the sense that one digs into other virtually hidden determinantal ideals coming from submatrices of the degenerate matrix.
These ideals come naturally while trying to uncover the nature of the relationship between the three ideals $J, I, J:I$.
One of the difficulties is that $I$ is not anymore prime for all values of $r$.
We conjecture that the bound ${{r+1}\choose {2}}\leq m-3$ is the exact obstruction for the primeness of $I$ (one direction is proved here).
The second conjectured statement is that the ring $R/J$ is Cohen--Macaulay if and only if $r=m-2$ (the ``only if'' part is proved here).

For a more precise discussion we refer to the statements of the various theorems.
As a guide, the main results are contained in Theorem~\ref{cloning_generic}, Theorem~\ref{primality_generic_cloned}, Theorem~\ref{dim_dual}, Theorem~\ref{expected_mult}, Theorem~\ref{polar-zeros}, Theorem~\ref{primality_generic_zeros} and Theorem~\ref{dual_zeros}.

Unless otherwise stated, we assume throughout that the ground field has characteristic zero.

\section{Preliminaries}

The aim of this section is to  review some notions and tools from ideal theory and its role in birational maps, including homaloidal ones.

\subsection{Review of ideal invariants}

Let $(R,\mathfrak{m})$ denote a Notherian local ring and its maximal ideal (respectively, a standard graded ring over a field and its irrelevant ideal).
For an ideal $I\subset \mathfrak{m}$ (respectively, a homogeneous ideal $I\subset \mathfrak{m}$), the \emph{special fiber} of $I$ is the ring $\mathcal{R}(I)/\mathfrak{m}\mathcal{R}(I)$.
Note that this is an algebra over the residue field of $R$.
The (Krull) dimension of this algebra is called the \emph{analytic spread} of $I$ and is denoted $\ell (I)$. 

Quite generally, given ideals $J\subset I$  in a ring $R$,  $J$ is said to be a \emph{reduction} of $I$ if there exists an integer $n\geq 0$ such that $I^{n+1}=JI^n.$
An ideal shares the same radical with all its reductions.
Therefore, they share the same set of minimal primes and have the same codimension.
A reduction $J$ of $I$ is called \emph{minimal} if no ideal strictly contained in $J$ is a reduction of $I$.
The \emph{reduction number} of $I$ with respect to a reduction $J$ is the minimum integer $n$ such that $JI^{n}=I^{n+1}$. It is denoted by $\mathrm{red}_{J}(I)$. The (absolute) \emph{reduction number} of $I$ is defined as $\mathrm{red}(I)=\mathrm{min}\{\mathrm{red}_{J}(I)~|~J\subset I~\mathrm{is}~\mathrm{a}~\mathrm{minimal}~\mathrm{reduction}~\mathrm{of}~I\}.$ 
If $R/\mathfrak{m}$ is infinite, then every minimal reduction of $I$ is minimally generated by exactly $\ell(I)$ elements. In particular, in this case, every reduction of $I$ contains a reduction generated by $\ell(I)$ elements.

The following invariants are related in the case of $(R,\mathfrak{m})$:
$$\mathrm{ht}(I)\leq \ell(I) \leq \min\{\mu(I), \mathrm{dim}(R)\},$$
where $\mu(I)$ stands for the minimal number of generators of $I$.
If the rightmost inequality turns out to be an equality, one says that $I$ has maximal analytic spread.
By and large, the ideals considered in this work will have $\dim R\leq \mu(I)$, hence being of maximal analytic spread means in this case that $\ell(I)=\dim R$.

Suppose now that  $R$ is a standard graded over a field $k$ and $I$ is minimally generated by $n+1$ forms of same degree $s$. In this case, $I$ is more precisely  given by means of a free graded presentation $$ R(-(s+1))^\ell \; \oplus\; \sum_{j\geq 2} R(-(s+j))\stackrel{\varphi}{\longrightarrow} R(-s)^{n+1}\longrightarrow I\longrightarrow 0  $$
for suitable shifts.   Of much interest in this work is the value of $\ell$.  The image of $R(-(s+1))^\ell$ by $\varphi$ is the {\it  linear part of $\varphi$} -- often denoted $\varphi_1$.  It is easy to see that the rank of $\varphi_1$  does not depend  on the particular minimal system of generators of $I$. Thus, we call it the {\it linear rank of} $I$. One says that $I$ has {\it maximal liner  rank} provided its linear rank is $n$ (=rank($\varphi$)). Clearly, the latter condition  is trivially satisfied if $\varphi=\varphi_1 $, in which case $I$ is said  to have {\it linear presentation} (or is {\it linearly presented}). 

Note that $\varphi$  is a graded matrix whose columns generate the (first) {\it syzygy module of $I$} (corresponding to the given choice of generators)  and  a {\it syzyzy} of $I$  is an element of this module -- that is, a linear  relation, with coefficients in $R$, on the chosen generators. In this context, $\varphi_1$ can be taken as the submatrix of $\varphi$  whose entries are linear forms of the standard graded ring $R$.  Thus, the linear rank is the rank of the matrix of the linear syzygies.

Recall the notion of the initial ideal of a polynomial ideal over a field.
For this one has to introduce a monomial order in the polynomial ring. Given such a monomial order, if $f\in R$ we denote by ${\rm in}(f)$ the initial term of $f$  and by ${\rm in}(I)$  the ideal generated by the initial terms of the elements of $I$ -- this ideal is called the initial ideal of $I$.

There are many excellent sources for the general theory of monomial ideals and Gr\"obner bases; we refer to the recent book \cite{HeHiBook}.

\subsection{Homaloidal polynomials}

Let $k$ be an arbitrary field. For the purpose of the full geometric picture we may assume $k$ to be algebraically closed. We denote by $\mathbb{P}^n=\mathbb{P}^n_k$ the $n$th projective space, where $n\geq 1$.

A rational map $\mathcal{F} : \mathbb{P}^n\dashrightarrow \mathbb{P}^m$  is defined by $m + 1$ forms ${\bf f} = \{f_0, \ldots , f_m\} \subset R :=k[{\bf x}] = k[x_0, \ldots, x_n]$ of the same degree $d\geq 1$, not all null.  We often write $\mathcal{F} = (f_0 : \cdots : f_m)$ to underscore the projective setup. Any rational map can without lost of generality be brought to satisfy the condition that ${\rm gcd}\{f_0,\ldots, f_m\} = 1$ (in the geometric terminology, $\mathcal{F}$ {\em has no fixed part}). The common degree $d$ of the $f_j$ is the degree  of $\mathcal{F}$ and the ideal $I_{\mathcal{F}}=\left(f_0,\ldots,f_m \right)$ is called the base ideal of $\mathcal{F}$.

The {\em image} of $\mathcal{F}$ is the projective subvariety $W \subset \mathbb{P}^m$ whose homogeneous coordinate ring is the $k$-subalgebra $k[{\bf f} ] \subset R$ after degree renormalization. Write $S := k[{\bf f} ] \simeq k[{\bf y}]/I(W)$, where $I(W) \subset k[{\bf y}] =k[y_0, \ldots , y_m]$ is the homogeneous defining ideal of the image in the embedding $W\subset \mathbb{P}^m$.  

We say that $\mathcal{F}$ is {\em birational onto its image} if there is a rational map $\mathcal{G}:\mathbb{P}^m\dashrightarrow \mathbb{P}^n$, say,  $\mathcal{G}=(g_0:\cdots :g_n)$, with the residue classes of the $g_i$'s modulo $I(W)$ not all vanishing, satisfying the relation
\begin{equation}
\label{composite_is_identity}
(g_0({\bf f}):\cdots : g_n({\bf f}))= (x_0:\cdots:x_n).
\end{equation}
(See \cite[Definition 2.10 and Corollary 2.12]{AHA}.)
When $m=n$ and $\mathcal{F}$  is a birational map of $\mathbb{P}^n$, we say that $\mathcal{F}$ is a {\em Cremona map}.  An important class of Cremona maps of $\mathbb{P}^n$ comes from the so-called  {\em polar maps}, that is, rational maps  whose coordinates are the partial derivatives  of a homogeneous polynomial $f$  in the ring $R=k[x_0,\ldots,x_n]$. More precisely:

Let $f\in k[{\bf x}]=k[x_0,\ldots,x_n]$ be a homogeneous polynomial of degree $d\geq 2$. The ideal 
	$$J=J_f=\left(\frac{\partial f}{\partial x_0},\ldots,\frac{\partial f}{\partial x_n}  \right)\subset k[{\bf x}]$$ 
is the Jacobian (or {\em gradient}) ideal of $f$. The rational map $\mathcal{P}_f:=\left( \frac{\partial f}{\partial x_0}:\cdots :\frac{\partial f}{\partial x_n}\right) $ is called  {\it the polar map} defined by $f$.  If  $\mathcal{P}_f$ is birational one says that $f$ is {\em homaloidal}.

We note  that the image of this map is the projective subvariety on the target whose homogeneous coordinate ring is given by the $k$-subalgebra $k[\partial f/\partial x_0,\ldots,\partial f/\partial x_n]\subset k[\xx]$ up to degree normalization. 
The image of $\mathcal{P}_f$ is called the {\it polar variety} of $f$.

The following  birationality criterion will be largely used in this work:

\begin{Theorem}{\rm[\cite{AHA}, Theorem 3.2]} \label{basic_criterion}
	Let $\mathcal{F}:\mathbb{P}^n\dashrightarrow \mathbb{P}^m$ be a rational map, given by $m+1$ forms ${\bf f} = \{f_0, \ldots , f_m\} $ of a fixed degree.  If $\dim (k[{\bf f}])=n+1$ and the linear rank of the base ideal $I_{\mathcal{F}}$  is $m$ {\rm ({\em maximal possible})} then $\mathcal{F}$  is birational onto its image.
\end{Theorem}

It is a classical result in characteristic zero that the Krull dimension of the $k$-algebra $k[{\bf f}]$  coincides with the rank of the Jacobian matrix of  ${\bf f} = \{f_0, \ldots , f_m\} $.  Assuming that the ground field has characteristic zero, the above criterion says that if the Hessian determinant $h(f)$ does not vanish  and the linear rank of the  gradient ideal of $f$ is maximal, then $f$ is homaloidal.

There are many sources for the basic material in these preliminaries; we refer to \cite{AHA}.

\section{Degeneration by cloning}\label{CLONING}

Quite generally, let $(a_{i,j})_{1\leq i\leq j\leq m}$ denote an $m\times m$ matrix where $a_{i,j}$ is either a variable on a ground polynomial ring $R=k[\xx]$ over a field $k$ or $a_{i,j}=0$.
One of the simplest specializations consists in going modulo a binomial of the shape $a_{i,j}-a_{i',j'}$, where $a_{i,j}\neq a_{i',j'}$ and  $a_{i',j'}\neq 0$.
The idea is to replace a certain nonzero entry $a_{i',j'}$ (variable)  by a different entry $a_{i,j}$, keeping $a_{i,j}$ as it was -- somewhat like {\em cloning} a variable and keeping the mold.
It  seems natural to expect that the new cloning place should matter as far as the finer properties of the ideals are concerned.

The main object of this section is the behavior of the generic square matrix under this sort of cloning degeneration.
We will use the following notation for the generic square matrix:
\begin{equation}\label{generic}
	\mathcal{G}:=\left(
	\begin{matrix}
		x_{1,1}&x_{1,2}&\ldots & x_{1,m-1} & x_{1,m}\\
		x_{2,1}&x_{2,2}&\ldots & x_{2,m-1} & x_{2,m}\\
		\vdots &\vdots &\ldots &\vdots & \vdots\\
		x_{m-1,1} & x_{m-1,2}&\ldots & x_{m-1,m-1} & x_{m-1,m}\\
		x_{m,1} & x_{m,2}&\ldots & x_{m,m-1} & x_{m,m}
	\end{matrix}
	\right),
\end{equation}
where the entries are independent variables over a field $k$.

Now, we distinguish essentially two sorts of cloning: the one that replaces an entry $x_{i',j'}$ by another entry $x_{i,j}$ such that $i\neq i'$ and $j\neq j'$, and the one in which this replacement has either $i= i'$ or $j= j'$.

In the situation of the second kind of cloning, by an obvious elementary operation and renaming of  variables  (which is possible since the original matrix is generic), one can assume that the matrix is the result of replacing a variable  by zero on a generic matrix.
Such a procedure is recurrent, letting several entries being replaced by zeros.
The resulting matrix along with its main properties will be studied in Section~\ref{zeros}.

Therefore, this section will deal exclusively with the first kind of cloning -- which, for emphasis, could be refereed to as {\em diagonal cloning}.
Up to elementary row/column operations and renaming of variables, we assume once for all that the diagonally cloned matrix has the shape
\begin{equation}\label{generic_cloned}
\mathcal{GC}:=\left(
\begin{matrix}
x_{1,1}&x_{1,2}&\ldots & x_{1,m-1} & x_{1,m}\\
x_{2,1}&x_{2,2}&\ldots & x_{2,m-1} & x_{2,m}\\
\vdots &\vdots &\ldots &\vdots & \vdots\\
x_{m-1,1} & x_{m-1,2}&\ldots & x_{m-1,m-1} & x_{m-1,m}\\
x_{m,1} & x_{m,2}&\ldots & x_{m,m-1} & x_{m-1,m-1}
\end{matrix}
\right),
\end{equation}
where the entry $x_{m-1,m-1}$ has been cloned as the $(m,m)$-entry of the $m\times m$ generic matrix.
The terminology may help us remind of the close interchange between properties associated to one or the other copy of the same variable in its place as an entry of the matrix. 
The question as to whether there is a similar theory for repeatedly many diagonal cloning steps has not been taken up in this work, but it looks challenging.

\smallskip

Throughout  $I_r({M})$ denotes the ideal generated by the $r$-minors  of a matrix ${M}$.

The following notion has been largely dealt with in \cite{Eisenbud2}.

An  $m\times n$  matrix ${M}$ of linear forms ($m\leq n$) over a ground field  is said to be {\em $s$-generic} for some integer $1\leq s\leq m$ if even after arbitrary invertible row and column operations, any $s$ of its entries are linearly independent over the field.	
It was proved in \cite{Eisenbud2} that the $m\times n$ generic matrix over a field  is $m$-generic; in particular, this matrix is $s$-generic for any $1\leq s\leq m$.  

Most specializations of the generic matrix fail to be $s$-generic for $s\geq 2$ due to their very format.
However many classical matrices are $1$-generic.

One of the important consequences of $s$-genericity is the  primeness of the ideal of $r$-minors for certain values of $r$. With an appropriate adaptation  of the original notation, the part of the result  we need  reads as follows:

\begin{Proposition}{\rm (\cite[Theorem 2.1]{Eisenbud2})}\label{Pisprime} One is given integers $1\leq w\leq v$.
	Let $\mathcal{G}$ denote the $w\times v$ generic matrix over a ground field. 
	Let ${M'}$ denote a $w\times v$ matrix of linear forms in the entries of $\mathcal{G}$ and let further ${M}$
	denote a $w\times v$ matrix of linear forms in the entries of ${M'}$.
	Let there be given an integer $k\geq 1$ such that ${M'}$ is a $(w-k)$-generic matrix and such that the vector space spanned by the entries of ${M}$ has codimension at most $k-1$ in the vector	space spanned by the entries of ${M'}$.
	Then the ideal $I_{k+1}(M)$ is  is prime.
\end{Proposition}

The following result originally appeared in \cite{Golberg} in a different context. It has independently been obtained in \cite[Proposition 5.3.1]{Maral} in the presently stated form.

\begin{Proposition}\label{GolMar}
	Let $M$  denote a square matrix over $R=k[x_0,\ldots,x_n]$ such that every entry is either $0$ or $x_i$ for some $i=1,\ldots,n.$ Then, for each $i=0,\ldots,n$, the partial derivative  of $f=\det M$ with respect to $x_i$ is the sum of the {\rm (}signed{\rm )} cofactors of the entry $x_i$, in all its slots as an entry of $M$.
\end{Proposition}

\subsection{Polar behavior}

Throughout we set $f:=\det (\mathcal{GC})$ and let $J=J_f\in R$ denote the gradient ideal of $f$, i.e., the ideal generated by the partial derivatives of $f$ with respect to the variables of $R$, the polynomial ring in the entries of $\mathcal{GC}$ over the ground field $k$.

\begin{Theorem}\label{cloning_generic}
Consider the diagonally cloned matrix as in {\rm (\ref{generic_cloned})}.
One has:
\begin{itemize}
\item[{\rm (i)}] $f$ is irreducible.
\item[{\rm (ii)}] The Hessian determinant $h(f)$ does not vanish.
\item[{\rm (iii)}] The linear rank of the gradient ideal of $f$ is $m^2-2$ {\rm (}maximum possible{\rm )}.
\item[{\rm (iv)}] $f$ is homaloidal.
\end{itemize}
\end{Theorem}
\demo
(i) We induct on $m$, the initial step of the induction being subsumed in the general step.

Expanding $f$ according to Laplace rule along the first row yields $$f=x_{1,1}\Delta_{1,1}+g,$$ where $\Delta_{1,1}$ is the determinant of the  $(m-1)\times (m-1)$  cloned generic matrix obtained from $\mathcal{GC}$ by omitting the first row  and the first column.
Note that both $\Delta_{1,1}$ and $g$ belong to the subring $k\left[ x_{1,2},\ldots,\ldots, x_{m,m-1}\right]$.
Thus, in order to show that $f$  is irreducible it suffices to prove that it is a primitive polynomial  (of degree 1) in $k\left[x_{1,2},\ldots ,x_{m,m-1} \right]\left[ x_{1,1}\right] $.

Now, on one hand, $\Delta_{1,1}$ is the determinant of a cloned matrix of the same type, hence it is irreducible by the inductive hypothesis. Therefore, it is enough to see that $\Delta_{1,1}$ is not factor of $g$. For this, one verifies their initial terms in the revlex monomial order, noting that they are slightly modified from the generic case: ${\rm in}(\Delta_{1,1})=(x_{2,m-1}x_{3,m-2}\cdots x_{m-1,2})x_{m-1,m-1}$ and ${\rm in}(g)={\rm in}(f)=(x_{1,m-1}x_{2,m-2}\cdots x_{m-1,1})x_{m-1,m-1}$.

An alternative more sophisticated argument is to use that the ideal $P$ of submaximal minors  has codimension $4$, as shown independently in Theorem~\ref{primality_generic_cloned} (i) below. 
Since $P=(J,\Delta_{m,m})$, as pointed out in the proof of the latter proposition, then $J$ has codimension at least $3$.
Therefore, the ring $R/(f)$ is locally regular in codimension one, so it must be normal.
But $f$ is homogeneous, hence irreducible.

\smallskip

(ii)
Set $\mathbf{v}:=\{x_{1,1},x_{2,2},x_{3,3},\ldots,x_{m-1,m-1}\}$ for the set of variables along the main diagonal.
We argue by a specialization procedure, namely, consider the ring endomorphism $\phi$ of $R$ by mapping any variable in  $\mathbf{v}$ to itself and by mapping any variable off  $\mathbf{v}$ to zero.
Clearly, it suffices to show that by applying $\phi$ to the entries of the Hessian matrix ${H}(f)$ the resulting matrix $\mathcal M$ has a nonzero determinant.

Note that the partial derivative of $f$ with respect to any  $x_{i,i}\in \mathbf{v}$ coincides with the signed cofactor of  $x_{i,i}$, for $i\leq m-2$, while for $i=m-1$ it is the sum of the respective signed cofactors of $x_{i,i}$ corresponding to its two appearances.

By expanding each such cofactor according to the Leibniz rule it is clear that it has a unique (nonzero) term whose support lies in $\mathbf{v}$ and, moreover, the remaining terms have degree at least $2$ in the variables off $\mathbf{v}$.

Now, for $x_{i,j}\notin \mathbf{v}$, without exception, the corresponding partial derivative coincides with the signed cofactor.
By a similar token, the Leibniz expansion of this cofactor  has no term whose support lies in $\mathbf{v}$ and has exactly one nonzero term of degree $1$ in the variables off $\mathbf{v}$.

By the preceding observation, applying $\phi$ to any second partial derivative of $f$ will return zero or a monomial supported on the variables in $\mathbf{v}$.
Thus, the entries of  $\mathcal M$ are zeros or monomials supported on the variables in $\mathbf{v}$.

To see that the determinant of the specialized matrix  $\mathcal M$ is nonzero, consider the Jacobian matrix of the set of partial derivatives $\{f_v\,|\, v\in\mathbf{v}\}$ with respect to the variables in $\mathbf{v}$.
Let $M_0$ denote the specialization of this Jacobian matrix by $\phi$, considered as a corresponding submatrix of  $\mathcal M$.
Up to permutation of rows and columns of  $\mathcal M$, we may write
$$\mathcal M=
\left(
  \begin{array}{cc}
    M_0 & N \\
    P & M_1
  \end{array}
\right),
$$
where $M_1$ has exactly one nonzero entry on each row and each column.
Now, by the way the second partial derivatives of $f$ specialize via $\phi$, as explained above, one must have $N=P=0$.
Therefore, $\det(\mathcal M)=\det(M_0)\det(M_1)$, so it remains to prove the nonvanishing of these two subdeterminants.

Now the first block $M_0$ is the Hessian matrix of the form
$$g:=\left(\prod_{i=1}^{m-2} x_{i,i}\right) x_{m-1,m-1}^2.
$$
This is the product of the generators of the $k$-subalgebra
$$k[x_{1,1},\ldots,x_{m-2,m-2}, x_{m-1,m-1}^2]\subset k[x_{1,1},\ldots,x_{m-2,m-2}, x_{m-1,m-1}].
$$
Clearly these generators are algebraically independent over $k$, hence the subalgebra is isomorphic to a polynomial ring itself.
Then $g$ becomes the product of the variables of a polynomial ring over $k$.
This is a classical homaloidal polynomial, hence we are done for the first matrix block.

As for $M_1$, since it has exactly one nonzero entry on each row and each column,
its determinant does not vanish.
\qed

\smallskip

(iii)  Let $f_{i,j}$ denote the  $x_{i,j}$-derivative of $f$ and let $\Delta_{j,i}$ stand for the (signed) cofactor of the $(i,j)$th entry of the matrix $\mathcal{GC}$.

The classical Cauchy cofactor formula
 \begin{equation}\label{CofactorFormula}
 \mathcal{GC} \cdot {\rm adj}(\mathcal{GC})={\rm adj}(\mathcal{GC})\cdot\mathcal{GC}=\det(\mathcal{GC})\, \mathbb{I}_m
 \end{equation}
yields by expansion a set of linear  relations involving the  (signed) cofactors of  $\mathcal{GC}$:
\begin{equation} \label{REL1}
\sum_{j=1}^m x_{i,j}\Delta_{j,k}=0,\ \mbox{for} \; 1\leq i \leq m-1 \; \mbox{and} \; 1\leq k  \leq m-2\;( k\neq i)
\end{equation}
\begin{equation}\label{REL2}
\sum_{j=1}^{m-1} x_{m,j}\Delta_{j,k}+x_{m-1,m-1}\Delta_{m,k}=0,\ \mbox{for} \;  1\leq k  \leq m-2
\end{equation}
\begin{equation}\label{REL3}
\sum_{j=1}^m x_{i,j}\Delta_{j,i}=\sum_{j=1}^m x_{i+1,j}\Delta_{j,i+1}, \; \mbox{for} \; 1\leq i\leq m-3
\end{equation}
\begin{equation}\label{REL4}
\sum_{i=1}^m x_{i,k}\Delta_{j,i}=0, \; \mbox{for} \; 1\leq j\leq m-3 \; \mbox{and} \; j< k \leq j+2.
\end{equation}
\begin{equation}\label{REL5}
\sum_{i=1}^m x_{i,m-1}\Delta_{m-2,i}=0. 
\end{equation}
\begin{equation}\label{REL6}
\sum_{i=1}^{m-1} x_{i,m}\Delta_{m-2,i}+x_{m-1,m-1}\Delta_{m-2,m}=0.
\end{equation}

Since $f_{i,j}=\Delta_{j,i}$  for every  $(i,j)\neq (m-1,m-1)$ and the above relations do not involve $\Delta_{m-1,m-1}$  or  $\Delta_{m,m}$ then they give linear syzygies of the partial derivatives of $f$.

In addition, (\ref{CofactorFormula}) yields the following linear relations:

\begin{equation}\label{rel1}
\sum_{j=1}^{m-1} x_{m-1,j} \Delta_{j,m} +x_{m-1,m} \Delta_{m,m}=0
\end{equation}
\begin{equation}\label{rel2}
\sum_{i=1}^{m-2} x_{i,m}\Delta_{m-1,i}+x_{m-1,m} \Delta_{m-1,m-1}+x_{m-1,m-1}\Delta_{m-1,m}=0
\end{equation}
\begin{equation}\label{rel3}
\sum_{i=1}^{m-1} x_{i,m-1}\Delta_{m,i}+x_{m,m-1}\Delta_{m,m} =0
\end{equation}
\begin{equation}\label{rel4}
\sum_{j=1}^{m-2}x_{m,j}\Delta_{j,m-1}+x_{m,m-1}\Delta_{m-1,m-1} +x_{m-1,m-1}\Delta_{m,m-1}=0
\end{equation}
\begin{equation}\label{rel5}
\sum_{j=1, {j\neq m-1}}^m x_{m-1,j}\Delta_{j,m-1}+x_{m-1,m-1}\Delta_{m-1,m-1}=\sum_{j=1}^m x_{m-2, j}\Delta_{j,m-2}
\end{equation}
\begin{equation}\label{rel6}
\sum_{j=1}^{m-1} x_{m,j}\Delta_{j,m}+x_{m-1,m-1}\Delta_{m,m}=\sum_{j=1}^m x_{m-2, j}\Delta_{j,m-2}.
\end{equation}
As  $f_{m-1, m-1}=\Delta_{m-1,m-1}+\Delta_{m,m}$, adding (\ref{rel1}) to (\ref{rel2}), (\ref{rel3}) to (\ref{rel4}) and (\ref{rel5}) to (\ref{rel6}), respectively, outputs three new linear syzygies of the partial derivatives of $f$.
Thus one has a total of $(m-2)(m-1)+(m-3)+2(m-2)+3=m^2-2$ linear syzygies of $J$.

It remains to show that these are independent.

For this we order the set of partial derivatives $f_{i,j}$ in accordance with the following ordered list of the entries $x_{i,j}$:
\begin{eqnarray} \nonumber
x_{1,1},x_{1,2},\ldots, x_{1,m}\rightsquigarrow x_{2,1},x_{2,2},\ldots, x_{2,m}\rightsquigarrow\ldots\rightsquigarrow x_{m-2,1},x_{m-2,2}\ldots, x_{m-2,m},\\ \nonumber 
x_{m-1,1},x_{m,1} \rightsquigarrow x_{m-1,2},x_{m,2}\rightsquigarrow\ldots\rightsquigarrow x_{m-1,m-1},x_{m,m-1}\rightsquigarrow x_{m-1,m}
\end{eqnarray}
Here we traverse the entries along the matrix rows, left to right, starting with the first row and stopping prior to the row having $x_{m-1,m-1}$ as an entry; then start traversing the last two rows along its columns top to bottom, until exhausting all variables.

We now claim that, ordering the set of partial derivatives $f_{i,j}$ in this way, the above sets of linear relations can be grouped into the following block matrix of  linear syzygies:

{\footnotesize
	$$\left(\begin{array}{cccc|ccccc|cccc}
	\varphi_1 &    &         &              &                 &     &        &                       &              &            &           &\\
	\mathbf{0}       &\varphi_2 &\ldots   &              &                 &     &        &                       &               &            &        &\\
	\vdots    &\vdots    & \ddots  & \vdots       &                 &     &        &                       &               &            &        &\\
	\mathbf{0}        & \mathbf{0}        & \ldots  & \varphi_{m-2}&                 &    &        &                       &      &   &        &\\ [2pt]
	\hline\\ [-9pt]
	\mathbf{0}_2^{m-1} & \mathbf{0}_2^m  & \ldots  &  \mathbf{0}_2^m	      & \varphi_{2}^3   &     &        &                       &               &            &         &\\
	\mathbf{0}_2^{m-1} & \mathbf{0}_2^m  &  \ldots & \mathbf{0}_2^m        &    \mathbf{0}_2^2         & \varphi_{3}^4   &  &                       &               &            &         &\\
	\vdots    & \vdots & \ldots  & \vdots       &  \vdots          &  \vdots        &\ddots & & & & &\\
	\mathbf{0}_2^{m-1} & \mathbf{0}_2^{m}& \ldots  &  \mathbf{0}_2^m	      &    \mathbf{0}_2^2        &  \mathbf{0}_2^2   & \ldots &   \varphi_{m-2}^{m-1}&            &            &        &\\
	\mathbf{0}_2^{m-1} & \mathbf{0}_2^m  & \ldots  &  \mathbf{0}_2^m	      &    \mathbf{0}_2^2        &  \mathbf{0}_2^2   & \ldots &   \mathbf{0}_2^2                  & \varphi_{m-1}^m &     &         &\\ [2pt]
	\hline\\ [-9pt]
	\mathbf{0}_1^{m-1}    & \mathbf{0}_1^{m}        & \ldots  &     \mathbf{0}_1^m  &    \mathbf{0}_1^2                 & \mathbf{0}_1^2      &    \ldots   &  \mathbf{0}_1^2        & \mathbf{0}_1^2     & x_{m-1,m} &  x_{m,m-1}   & x_{m-1,m-1} \\
	\mathbf{0}_ 1^{m-1}   & \mathbf{0}_1^m        & \ldots  &     \mathbf{0}_1^m    &    \mathbf{0}_1^2                & \mathbf{0}_1^2      &          \ldots   &  \mathbf{0}_1^2       &    \mathbf{0}_1^2    & 2x_{m-1,m-1} &   \mathbf{0} & x_{m,m-1} \\
	\mathbf{0}_1^{m-1}    & \mathbf{0}_1^m       & \ldots  & \mathbf{0}_1^m         & \mathbf{0}_1^2                   & \mathbf{0}_1^2     &         \ldots    &  \mathbf{0}_1^2     &  \mathbf{0}_1^2           & \mathbf{0}   &   2x_{m-1,m-1} & x_{m-1,m}\\
	\end{array}
	\right).$$
}
Let us explain the blocks of the above matrix:
\begin{itemize}
	\item $\varphi_1$ is the matrix obtained  from the  transpose $\mathcal{GC}^t$ of $\mathcal{GC}$ by omitting the first column;
	\item $\varphi_2,\ldots,\varphi_{m-2}$ are each a copy of  $\mathcal{GC}^t$ (up to column permutation);
	\item {\small$\varphi_{r}^{r+1} =\left(
		\begin{matrix}
		x_{m-1,r} & x_{m-1,r+1}\\
		x_{m,r}& x_{m,r+1)}\\
		\end{matrix}
		\right), $} $r=2,\ldots, m-2;$ {\small$
		\varphi_{m-1}^m=\left(
		\begin{matrix}
		x_{m-1,m-1} & x_{m-1,m}\\
		x_{m,m-1} & x_{m-1,m-1}
		\end{matrix}
		\right)$};
	\item Each $\mathbf{0}$ under $\varphi_1$  is an $m\times (m-1)$ block of zeros  and each $\mathbf{0}$ under $\varphi_i$ is an $m\times m$ block of zeros, for $i=2,\ldots,m-3$ ;
	\item $\mathbf{0}_r^c$ denotes  an $r\times c$ block of zeros, for  $r=1,2$ and $c=2,m-1, m$.
\end{itemize}

Next we justify why these blocks make up (linear) syzygies.

First, as already observed,  the relations (\ref{REL1}) through (\ref{rel6}) yield linear syzygies of the partial derivatives of $f$.
Setting $k=1$ in the  relations  (\ref{REL1}) and (\ref{REL2}) the resulting expressions can be written, respectively,  as \begin{center}$\sum_{j=1}^m x_{i,j}f_{1,j}=0$ for all  $i=2,\ldots, m-1$ and $\sum_{j=1}^{m-1} x_{m,j}f_{1,j}+x_{m-1,m-1}f_{1,m}=0$.\end{center}  
Ordering the set of partial derivatives $f_{i,j}$ as explained before, the coefficients of these relations form the first matrix above
$$\varphi_1:=\left(\begin{array}{ccccc} 
	x_{2,1} & x_{3,1}& \ldots  &x_{m-1, 1} & x_{m,1}\\
	\vdots  & \vdots & \ldots  &\vdots  & \vdots \\
	x_{2,m-1} & x_{3,m-1}& \ldots  &x_{m-1, m-1} & x_{m,m-1}\\
	x_{2,m} & x_{3,m}& \ldots  &x_{m-1, m} & x_{m-1,m-1}
	\end{array}
	\right)
	$$
Note that $\varphi_1$ coincides indeed with the submatrix of $\mathcal{GC}^t$ obtained by omitting its first column. 
	
Getting $\varphi_k$, for  $k=2,\ldots.,m-2$, is similar, namely, use again relations  (\ref{REL1}) and (\ref{REL2}) retrieving a submatrix of $\mathcal{GC}^t$ excluding the $k$th column and replacing it with an extra column that comes from relation (\ref{REL3}) taking $i=k-1$.

Continuing, for each $r=2,\ldots ,m-2$ the  block $\varphi_{r}^{r+1}$ comes from the relation (\ref{REL4}) (setting $j=r-1$) and $\varphi_{m-1}^{m}$ comes from  the relations (\ref{REL5}) and (\ref{REL6}). Finally, the lower right corner $3\times 3$ block of the matrix of  linear syzygies comes from the three last relations 
obtained by adding (\ref{rel1}) to (\ref{rel2}), (\ref{rel3}) to (\ref{rel4}) and (\ref{rel5}) to (\ref{rel6}).

This proves the claim about the large matrix above.
Counting through the sizes of the various blocks, one sees that this matrix is $(m^2-1)\times (m^2-2)$.
Omitting its first row obtains a block-diagonal submatrix of size
$(m^2-2)\times (m^2-2)$, where each block has nonzero determinant. Thus, the linear rank of $J$ attains the maximum.

\smallskip

(iv) By (ii) the polar map of $f$ is dominant. Since the linear rank is maximum by (iii), one can apply Theorem~\ref{basic_criterion} to conclude that $f$ is homaloidal.
\qed

\medskip

\subsection{The ideal of the submaximal minors}

In this part we study the nature of the ideal of submaximal minors (cofactors) of $\mathcal{GC}$.
As previously, $J$ denotes the gradient ideal of $f=\det (\mathcal{GC})$

\begin{Theorem}\label{primality_generic_cloned} 
Consider the matrix $\mathcal{GC}$ as in {\rm (\ref{generic_cloned})}, with $m\geq 3$.
Let $P:=I_{m-1}(\mathcal{GC})$ denote the ideal of $(m-1)$-minors of $\mathcal{GC}$.
Then
\begin{enumerate}
\item[{\rm (i)}]  $P$ is a Gorenstein prime ideal of codimension $4$.
\item[{\rm (ii)}]  $J$ has codimension $4$ and $P$ is the minimal primary component of $J$ in $R$. 
\item[{\rm (iii)}] $J$ defines a double structure on the variety defined by $P$, admitting one single embedded component, the latter being a linear space of codimension $4m-5$.
\item[{\rm (iv)}]  Letting $\mathbb{D}_{i,j}$ denote the cofactor of the $(i,j)$-entry of the generic matrix $(y_{i,j})_{1\leq i,j\leq m}$,
the $(m-1)$-minors $\boldsymbol\Delta=\{\Delta_{i,j}\}$ of $\mathcal{GC}$ define a birational map $\pp^{m^2-2}\dasharrow \pp^{m^2-1}$ onto a hypersurface of degree $m-1$ with defining equation $\mathbb{D}_{m,m}-\mathbb{D}_{m-1,m-1}$ and inverse map defined by the linear system spanned by  $\widetilde{\mathbb{D}}:=\{\mathbb{D}_{i,j}\,|\, (i,j)\neq (m,m)\}$ modulo $\mathbb{D}_{m,m}-\mathbb{D}_{m-1,m-1}$. 
\item[{\rm (v)}]  $J$ is not a reduction of $P$.
\end{enumerate}
\end{Theorem}
\demo
(i)  Let $\mathcal{P}$ denote the ideal of submaximal minors of the fully generic matrix (\ref{generic}). The linear form $x_{m,m}-x_{m-1,m-1}$ is regular on the corresponding polynomial ambient and also modulo $\mathcal{P}$ as the latter is prime and generated in degree $m-1\geq 2$.
Since $\mathcal{P}$ is a Gorenstein ideal of codimenson $4$ by a well-known result (``Scandinavian complex''), then so is $P$. 

In order to prove primality, we first consider the case $m=3$ which seems to require a direct intervention.
We will show more, namely, that $R/P$ is normal -- and, hence a domain as $P$ is a homogeneous ideal.
Since $R/P$ is a Gorenstein ring, it suffices to show that $R/P$ is locally regular in codimension one. For this consider the Jacobian matrix of $P$:

$$\left( \begin{matrix}
	x_{2,2} & -x_{2,1} & 0 & -x_{1,2} & x_{1,1} & 0 & 0 & 0\\
	x_{2,3} & 0 & -x_{2,1} & -x_{1,3} & 0 & x_{1,1} & 0 & 0\\
	0 & x_{2,3} & -x_{2,2} & 0 & -x_{1,3} & x_{1,2} & 0 & 0\\
	x_{3,2} & -x_{3,1} & 0 & 0 & 0 & 0 & -x_{1,2} & x_{1,1}\\
    x_{2,2} & 0 & -x_{3,1} & 0 & x_{1,1} & 0 & -x_{1,3} & 0\\
    0 & x_{2,2} & -x_{3,2} & 0 & x_{1,2} & 0 & 0  & -x_{1,3}\\
    0 & 0 & 0 & x_{3,2} & -x_{3,1} & 0 & -x_{2,2} & x_{2,1}\\
    0 & 0 & 0 & x_{2,2} & x_{2,1} & -x_{3,1} & -x_{2,3} & 0\\
    0 & 0 & 0 & 0 & 2x_{2,2} & -x_{3,2} & 0 & -x_{2,3}\\
\end{matrix} \right).$$

Direct inspection yields that the following pure powers are (up to sign) $4$-minors of this matrix: $x_{1,3}^4$, $x_{2,1}^4 $, $x_{2,2}^4$, $x_{2,3}^4$, $x_{3,1}^4$ and $x_{3,2}^4$. Therefore, the ideal of $4$-minors of the Jacobian matrix has codimension at least $6=4+2$, thus  ensuring that  $R/P$  satisfies $(R_1)$. 

\smallskip

For $m\geq 4$ we apply Proposition~\ref{Pisprime} with $M'=\mathcal{G}$ standing for an $m\times m$ generic matrix and $M=\mathcal{GC}$ the cloned generic matrix as in the statement.
In addition, we take $k=m-2$, so $k+1=m-1$ is the size of the submaximal minors.
Since $m\geq 4$ and the vector space codimension in the theorem is now $1$, one has $1\leq m-3=k-1$ as required.
Finally, the $m\times m$ generic matrix is $m$-generic as explained in \cite[Examples, p. 548]{Eisenbud2}; in particular, it is $2=m-(m-2)$-generic.
The theorem applies to give that the ideal $P=I_{m-1}(\mathcal{GC})$ is  prime.

\medskip

(ii) 
By item (i), $P$ is a prime ideal of codimension $4$.
We first show that ${\rm cod}(J:P)>4$, which ensures that the radical of the unmixed part of $J$ has no primes of codimension $<4$ and coincides with $P$ -- in particular, $J$ will turn out to have codimension $4$ as stated.

For this note that $P=(J,\Delta_{m,m})$, where $\Delta_{m,m}$ denotes the cofactor of the $(m,m)$th entry. From the cofactor identity we read the following relations:
\begin{eqnarray} \nonumber
\sum_{j=1}^m x_{k,j}\Delta_{j,m}&=&0,\; \mbox{for} \;  k=1,\ldots, m-1;\\\nonumber
\sum_{j=1}^{m-1} x_{m,j}\Delta_{m,j}+x_{m-1,m-1}\Delta_{m,m}&=&\sum_{j=1}^m x_{1,j}\Delta_{j,1}; \\ \nonumber
\sum_{i=1}^m x_{i,k}\Delta_{m,i}&=&0, \; \mbox{for} \;  k = 1,\ldots,m-1 ;
\end{eqnarray}

Since the partial derivative $f_{i,j}$ of $f$ with respect to the variable $x_{i,j}$ is the (signed) cofactor $\Delta_{j,i}$, with the single exception of the partial derivative with respect to the variable $x_{m-1,m-1}$,  we have that the entries of the $m$th column and of the $m$th row all belong to the ideal $\left(J:\Delta_{m,m} \right)=\left(J:P \right)$. In particular, the codimension of $\left(J:P \right) $ is at least $5$, as needed.

In addition, since $P$ has codimension $4$ then $J:P\not\subset P$.
Picking an element $a\in J:P\setminus P$ shows that $P_P\subset J_P$.
 Therefore $P$ is the unmixed part of $J$.

To prove that $P$ is actually the entire minimal primary component of $J$ we argue as follows.
In addition, also note that $P=(J,\Delta_{m-1,m-1})$,  where $\Delta_{m-1,m-1}$ denotes the cofactor of the  $(m-1,m-1)$th entry. From the cofactor identity we read the following relations:

\begin{eqnarray} \nonumber
\sum_{j=1, j\neq m-1}^m x_{k,j}\Delta_{j,m-1}+ x_{k,m-1}\Delta_{m-1,m-1}&=&0,\; \mbox{for} \;  k=1,\ldots, m, \;( k\neq m-1);\\\nonumber
\sum_{j=1, j\neq m-1}^m x_{m-1,j}\Delta_{j,m-1}+x_{m-1,m-1}\Delta_{m-1,m-1}&=&\sum_{j=1}^m x_{1,j}\Delta_{j,1}; \\ \nonumber
\sum_{i=1,i\neq m-1}^m x_{i,k}\Delta_{m-1,i}+x_{m-1,k}\Delta_{m-1,m-1}&=&0, \; \mbox{for} \;  k = 1,\ldots,m \;(k\neq m-1) ;
\end{eqnarray}
Then as above we have that the entries of the $(m-1)$th column and the $(m-1)$th row belong to the ideal $\left(J:\Delta_{m-1,m-1} \right)=\left(J:P \right)$. 

From this, the variables of the last two rows and columns of $\mathcal{GC}$ multiply $P$ into $J$. As is clear that $P$ is contained in  the ideal generated by these variables it follows that $P^2\subset J$ (of course, this much could eventually be verified by inspection). Therefore, the radical of $J$ -- i.e., the radical of the minimal primary part of $J$ -- is $P$.

(iii)  By (ii), $P$ is the minimal component of a primary decomposition of $J$.
 We claim that $J:P$ is generated by the $4m-5$ entries of $\mathcal{GC}$ off the upper left submatrix of size $(m-2)\times (m-2)$.
 Let $I$ denote the ideal generated by these entries.
 
 As seen in the previous item,  $I\subset J:P$. We now prove the reverse inclusion
by writing $I=I'+I''$ as sum of two prime ideals, where  $I'$ (respectively, $I''$) is the ideal generated by the variables on the $(m-1)$th row and on the $(m-1)$th column of $\mathcal{GC}$ (respectively, by the variables on the $m$th row and on the $m$th column of $\mathcal{GC}$). 
Observe that the cofactors $\Delta_{i,j}\in I''$ for all $(i,j)\neq (m,m)$ and  $\Delta_{i,j}\in I'$ for all $(i,j)\neq (m-1,m-1)$. Clearly, then $\Delta_{m,m}\notin I''$ and $\Delta_{m-1,m-1}\notin I'$.

Let $b\in J:P=J:  \Delta_{m,m}$, say,
  \begin{eqnarray} \nonumber b\,\Delta_{m,m} &=&\sum_{(i,j)\neq(m-1,m-1)} a_{i,j}f_{i,j}+af_{m-1,m-1}\\ \label{Q1}
  &=&\sum_{(i,j)\neq(m-1,m-1)} a_{i,j}\Delta_{j,i}+a(\Delta_{m-1,m-1}+\Delta_{m,m})
  \end{eqnarray}
for certain $a_{i,j},a\in R$. 
Then 
$$ (b-a)\Delta_{m,m} =\sum_{(i,j)\neq(m-1,m-1)} a_{i,j}\Delta_{j,i}+a\Delta_{m-1,m-1}  \in I''.$$ 
Since $I''$ is a prime ideal and $\Delta_{m,m}\notin I''$, we have $c:=b-a\in I''$.  
Substituting for $a=b-c$ in  (\ref{Q1}) gives   
$$  (-b+c)\Delta_{m-1,m-1}=\sum_{(i,j)\neq(m-1,m-1)} a_{i,j}\Delta_{j,i}-c\Delta_{m,m}\in I'.$$    
  
By a similar token, since $\Delta_{m-1,m-1}\notin I'$, then $-b+c\in I'$. Therefore $$b=c-(-b+c)\in I''+I'=I,$$
as required.
  
In particular, $J:P$ is a prime ideal which is necessarily an associated prime of prime of $R/J$. As pointed out, $P\subset J:P$, hence $J:P$ is an embedded prime of $R/J$. Moreover, this also gives $P^2\subset J$, hence $J$ defines a double structure on the irreducible variety defined by $P$.

Let $\mathcal{Q}$ denotes the embedded component of $J$ with radical $J:P$ and let 
$\mathcal{Q}'$ denote the intersection of the remaining embedded components of $J$.
From $J=P\cap \mathcal{Q}\cap \mathcal{Q}'$ we get 
$$J:P=(\mathcal{Q}:P)\cap (\mathcal{Q}':P),$$
in particular, passing to radicals,  $J:P\subset \sqrt{\mathcal{Q}'}$.
This shows that $\mathcal{Q}$ is the unique embedded component of codimension $\leq 4m-5$ and the corresponding geometric component is supported on a linear subspace.

\medskip

(iv) By Theorem~\ref{cloning_generic}~(ii), the polar map is dominant, i.e., the partial derivatives of $f$ generate a subalgebra of maximum dimension  ($=m^2-1$).
Since $J\subset P$ is an inclusion in the same degree, the subalgebra generated by the submaximal minors has dimension $m^2-1$ as well.
On the other hand, since $P$ is a specialization from the generic case, it is linearly presented.
Therefore, by Theorem~\ref{basic_criterion} the minors define a birational map onto a hypersurface.

To get the inverse map and the defining equation of the image we proceed as follows.

\smallskip

Write $\Delta_{j,i}$ for the cofactor of the $(i,j)$-entry of $\mathcal{GC}$.
For the image it suffices to show that $\mathbb{D}_{m,m}-\mathbb{D}_{m-1,m-1}$ belongs to the kernel of the $k$-algebra map 
$$\psi:k[y_{i,j}\,|\, 1\leq i,j\leq m]\rar k[\underline{\Delta}]=k[\Delta_{i,j}\,|\, 1\leq i,j\leq m],$$ 
as it is clearly an irreducible polynomial.

Consider the following well-known matrix identity
	\begin{equation}\label{inversion_formula}
	{\rm adj}({\rm adj}(\mathcal{GC}))=f^{m-2}\cdot\mathcal{GC},
	\end{equation}
where adj$(M)$ denotes the transpose matrix of cofactors of a square matrix $M$.
On the right-hand side matrix we obviously see the same element as its $(m-1,m-1)$-entry as its $(m,m)$-entry, namely, $f^{m-2} x_{m-1,m-1}$.

As to the entries of the matrix on the left-hand side, for any $(k,l)$, the $(k,l)$-entry is $\mathbb{D}_{l,k}(\underline{\Delta})$.
Indeed, the $(k,l)$-entry of ${\rm adj}({\rm adj}(\mathcal{GC}))$ is the cofactor of the entry $\Delta_{l,k}$ in the matrix ${\rm adj}(\mathcal{GC})$. Clearly, this cofactor is the $(l,k)$-cofactor $\mathbb{D}_{l,k}$ of the generic matrix $(y_{i,j})_{1\leq i,j\leq m}$ evaluated at $\underline{\Delta}$. 

Therefore, we get $(\mathbb{D}_{m,m}-\mathbb{D}_{m-1,m-1})(\underline{\Delta})=0$, as required.

Finally, by the same token, from (\ref{inversion_formula}) one deduces that the inverse map has coordinates  $\widetilde{\mathbb{D}}:=\{\mathbb{D}_{i,j}\,|\, (i,j)\neq (m,m)\}$ modulo $\mathbb{D}_{m,m}-\mathbb{D}_{m-1,m-1}$.

\smallskip

(v) It follows from (iv) that the reduction number of a minimal reduction of $P$ is $m-2$.
Thus, to conclude, it suffices to prove that $P^{m-1}\not\subset JP^{m-2}$.

We will show that $\Delta^{m-1}\in P^{m-1}$ does not belong to  $JP^{m-2}$.

Recall from previous passages that $J$ is generated by the cofactors 
$$\Delta_{l,h},\; \mbox{\rm with} \;(l,h)\neq (m-1,m-1), (l,h)\neq (m,m)
$$
and the additional form $\Delta_{m,m}+\Delta_{m-1,m-1}$.

If  $\Delta^{m-1}\in  JP^{m-2}$, we can write

\begin{equation}\label{reduction-relation}
\Delta_{m,m}^{m-1}=\sum_{(l,h)\neq (m-1,m-1)\atop (l,h)\neq (m,m)} \Delta_{l,h}Q_{l,h}(\underline{\Delta})+(\Delta_{m,m}+\Delta_{m-1,m-1})Q(\underline{\Delta})
\end{equation}

\noindent where $Q_{l,h}(\underline{\Delta})$ and $Q(\underline{\Delta})$ are homogeneous polynomial expressions of degree $m-2$ in the set $$\underline{\Delta}=\{\Delta_{i,j}\,| \,1\leq i\leq j\leq m\}$$ of the cofactors (generators of $P$). 

Clearly, this gives a polynomial relation of degree $m-1$ on the generators of $P$, so the corresponding form of degree $m-1$ in $k[y_{i,j}| 1\leq i\leq j\leq m]$ is a scalar multiple of the defining equation $\boldsymbol H:=\mathbb{D}_{m,m}-\mathbb{D}_{m-1,m-1}$ obtained in the previous item.
Note that $\mathbf{H}$ contains only squarefree terms.
We now argue that such a relation is impossible.

Namely, observe that the sum $$\sum_{(l,h)\neq (m-1,m-1)\atop (l,h)\neq (m,m)} \Delta_{l,h}Q_{l,h}(\underline{\Delta})$$  does not contain any nonzero terms of the form $\alpha\Delta_{m,m}^{m-1}$ or $\beta\Delta_{m-1,m-1}\Delta_{m,m}^{m-2}$.
In addition, if these two terms appear in $(\Delta_{m,m}+\Delta_{m-1,m-1})Q(\underline{\Delta})$ they must have the same scalar coefficient, say, $c\in k$. 
Bring the first of these to the left-hand side of (\ref{reduction-relation}) to get a polynomial relation of $P$ having a nonzero term $(1-c)y_{m-1,m-1}^{m-1}$.
If $c\neq 1$, this is a contradiction due to the squarefree nature of $\mathbf{H}$.

On the other hand, if $c=1$ then we still have a polynomial relation of $P$ having a nonzero term $y_{m-1,m-1}y_{m,m}^{m-2}$.
Now, if $m>3$ this is again a contradiction vis-\`a-vis the nature of $\mathbf{H}$ as the nonzero terms of the latter are squarefree monomials of degree $m-1>3-1=2$.
Finally, if $m=3$ a direct checking shows that the
 monomial $y_{m-1,m-1}y_{m,m}$ cannot be the support of a nonzero term in $\boldsymbol H$. 
 This concludes the statement. 
\qed

\subsection{The dual variety}

An interesting question in general is whether $f$ is a factor of its Hessian determinant $h(f)$ with multiplicity $\geq 1$. 
If this is the case, then $f$ is said in addition to have the {\em expected multiplicity} (according to Segre) if its multiplicity as a factor of $h(f)$ is  $m^2-2-\dim V(f)^*-1=m^2-3-\dim V(f)^*= {\rm cod}(V(f)^*)-1$, where $V(f)^*$ denotes the dual variety to the hypersurface $V(f)$ (see \cite{CRS}).

\begin{Theorem}\label{dim_dual} Let $f=\det(\mathcal{GC})$. 
	Then $\dim V(f)^*= 2m-2$. In particular, the expected multiplicity of $f$ as a factor of $h(f)$ is $m(m-2)-1$.
\end{Theorem}
\demo We develop the argument in two parts:

{\bf 1}. $\dim V(f)^*\geq  2m-2$. 

We draw on a result of Segre (\cite{Segre}), as transcribed in \cite[Lemma 7.2.7]{Frusso}, to wit:
$$\dim V(f)^*=\rk H(f) \pmod{f}-2,$$ 
where $H(f)$ denotes the Hessian matrix of $f$.
It will then suffice to show that $H(f)$ has a submatrix of rank at least $2m$ modulo $f$.
Consider the submatrix
$$\phi=\left(\begin{array}{ccccccc}
x_{2,1}&\ldots& x_{2,m-1}&x_{2,m}\\
\vdots&\ddots& \vdots&\vdots\\
x_{m-1,1}&\ldots&x_{m-1,m-1}&x_{m-1,m}\\
x_{m,1}&\ldots&x_{m,m-1}&x_{m-1,m-1}
\end{array}\right)$$
of $\mathcal{GC}$ obtained by omitting the first row. 
The maximal minors of this $(m-1)\times m$ matrix generate a codimension $2$ ideal.
On the other hand, $\phi$ has the property $F_1$ for its Fitting ideals.
Indeed, the Fitting ideals of its generic predecessor are prime ideals, hence the Fitting ideals of $\phi$ are specializations thereof, and as such each has the same codimension as the respective predecessor.
It follows from this that the ideal is of linear type (\cite{Trento}).
In particular the $m$ maximal minors of $\phi$ are algebraically independent over $k$, hence their Jacobian matrix with respect to the entries of $\phi$ has rank $m$.

Let $A'$ denote an $m\times m$ submatrix thereof with $\det(A')\neq 0$.
Now, $f\in(x_{1,1},\ldots,x_{1,m})$ while $\det(A')\notin (x_{1,1},\ldots,x_{1,m}).$
This means that $\det(A')\neq 0$ even modulo $f$.

Write the Hessian matrix $H(f)$ in the block form
$$H(f)=\left(\begin{array}{cccc}
\boldsymbol0&A\\
A^t&B
\end{array}\right),$$
where the first block row is the Jacobian matrix of the maximal minors of $\phi$ in the order of the variables starting with $\{x_{1,1},\ldots,x_{1,m}\}$ and $A^t$ denotes the transpose of $A$.
Since $A'$ above is a submatrix of $A$  of rank $m$ modulo $f$, then $h(f)=\det H(f)$ has rank at least $2m$ modulo $f$.

\smallskip

{\bf 2.} $\dim V(f)^*\leq  2m-2$.

Here we focus on the homogeneous coordinate ring of $V(f)^*$, namely, the following $k$-subalgebra of $k[\xx]/(f)$
\begin{equation}\label{algevra_of_dual}
k[\partial f/\partial x_{1,1},\ldots, \partial f/\partial x_{m,m-1}]/(f)\simeq k[y_{i,j}\,|\, 1\leq i\leq j\leq m, (i,j)\neq (m,m)]/P,
\end{equation}
for a suitable prime ideal $P$, the homogeneous defining ideal of $V(f)^*$.
The isomorphism is an isomorphism of graded $k$-algebras induced by the assignment $y_{i,j}\mapsto \partial f/\partial x_{i,j},\, (i,j)\neq (m,m)$. 

\smallskip

{\sc Claim 1.} The homogeneous defining ideal $P$ of $V(f)^*$ contains the ladder determinantal ideal generated by the $2\times 2$ minors of the following matrix
$$\mathcal{L}=\left(\begin{array}{cccccccc}
y_{1,1}&y_{1,2}&\ldots& y_{1,m-2}&y_{1,m}&y_{1,m-1}\\
y_{2,1}&y_{2,2}&\ldots& y_{2,m-2}&y_{2,m}&y_{2,m-1}\\
\vdots&\vdots&\ddots&\vdots&\vdots\\
y_{m-2,1}&y_{m-2,2}&\ldots& y_{m-2,m-2}&y_{m-2,m}&y_{m-2,m-1}\\
y_{m-1,1}&y_{m-1,2}&\ldots& y_{m-1,m-2}&y_{m-1,m}&\\
y_{m,1}&y_{m,2}&\ldots& y_{m,m-2}
\end{array}\right).$$

To see  this, we first recall that, for $(i,j)\neq (m-1,m-1)$, the partial derivative $\partial f/\partial x_{i,j}$ coincides with the cofactor of $x_{i,j}$ in $\mathcal{GC}$.
Since $y_{m-1,m-1}$ and $y_{m,m}$ are not entries of $\mathcal{L}$, the polynomial relations of the partial derivatives possibly involving the variables which are entries of $\mathcal{L}$ are exactly relations of the cofactors other than the cofactor of $x_{m-1,m-1}$.

Thus, we focus on these factors, considering the following relation afforded by the cofactor identity:
\begin{equation}\label{relation_mod_f}
{\rm adj}(\mathcal{GC})\cdot\mathcal{GC}\equiv 0\,({\rm mod\,}f).
\end{equation}
Further, for each pair of integers $i,j$ such that $1\leq i<j\leq m $ let $F_{ij}$ denote the $2\times m$ submatrix of ${\rm adj}(\mathcal{GC})$ consisting of the $i$th and $j$th rows.
In addition, let $C$ stand for the $m\times (m-1)$ submatrix of $\mathcal{GC}$  consisting of its $m-1$ leftmost columns.
Then (\ref{relation_mod_f}) gives the relations
$$F_{ij}C\equiv0\,({\rm mod\,}f),$$
for all $1\leq i<j\leq m $.
From this, since the rank of $C$ modulo $(f)$ is obviously still $m-1$, the one of every $F_{i,j}$ is necessarily $1$.
This shows that every $2\times 2$ minor of ${\rm adj}(\mathcal{GC})$ vanishes modulo $(f)$.
Therefore, every such minor that does not involve either one of the cofactors $\Delta_{m-1,m-1}$ and $\Delta_{m,m}$ gives a $2\times 2$ minor of $\mathcal{L}$ vanishing on the partial derivatives. Clearly, by construction, we obtain this way all the $2\times 2$ minors of $\mathcal{L}$.
This proves the claim.

\smallskip

Now, since $I_2(\mathcal{L})$ is a ladder determinantal ideal on a suitable  generic matrix it is a Cohen-Macaulay prime ideal (see \cite{Nar} for primeness and \cite{HeTr} for Cohen--Macaulayness). 
Moreover, its codimension is $m(m-2)-2=m^2-3-(2m-1)$ as follows from an application to this case of the general principle in terms of maximal chains as described in \cite[Theorem 4.6 and Corollary 4.7]{HeTr}.

Note that by $I_2(\mathcal{L})$ we understand the ideal generated by the $2\times 2$ minors of $\mathcal{L}$ in the polynomial ring $A:=k[I_1(\mathcal{L})_1]$ spanned by the entries of $\mathcal{L}$.
Clearly, its extension to the full polynomial ring $B:=k[y_{i,j}\,|\, 1\leq i\leq j\leq m, (i,j)\neq (m,m)]$, is still prime of codimension $m(m-2)-2$.
Thus, $P$ contains a prime subideal of codimension $m(m-2)-2$.

In addition, direct checking shows that the following two quadrics 
\begin{equation*}
g:=y_{1,1}y_{m,m-1}-y_{1,m-1}y_{m,1},\, h:=y_{1,1}y_{m-1,m-1}-y_{m-1,1}y_{1,m-1}-y_{m,1}y_{1,m}
\end{equation*}
belong to  $P$ and, furthermore:

{\sc Claim 2.} $\{g,h\}$ is a regular sequence modulo  $I_2(\mathcal{L})B$.

It suffices to prove the assertion locally at the powers of $y_{1,1}$ since $I_2(\mathcal{L})B$ is prime and $y_{1,1}\notin I_2(\mathcal{L})B$.
Now, locally at $y_{1,1}$ and at the level of the ambient rings, $f,g$ are like the variables $y_{m,m-1}$ and $y_{m-1,m-1}$, hence one has
$$B_{y_{_{1,1}}}/(g,h)\simeq A_{y_{_{1,1}}}.$$
Setting $I=I_2(\mathcal{L})$ for lighter reading, it follows that
\begin{eqnarray*}\nonumber
	(B_{y_{_{1,1}}}/IB_{y_{_{1,1}}})/(g,h)(B_{y_{_{1,1}}}/IB_{y_{_{1,1}}})&\simeq & B_{y_{_{1,1}}}/(g,h, I)\\ \nonumber
	&\simeq & (B_{y_{_{1,1}}}/(g,h))/((g,h, I)/(g,h))\\ \nonumber
	&\simeq & A_{y_{_{1,1}}}/IA_{y_{_{1,1}}}.
\end{eqnarray*}
But clearly, $B_{y_{_{1,1}}}/IB_{y_{_{1,1}}}\simeq (A_{y_{_{1,1}}}/IA_{y_{_{1,1}}})[ y_{m-1,m-1},y_{m,m-1}]_{y_{_{1,1}}}$, hence
$$\dim A_{y_{_{1,1}}}/IA_{y_{_{1,1}}}= \dim B_{y_{_{1,1}}}/IB_{y_{_{1,1}}}-2.$$
This shows that 
$$\dim (B_{y_{_{1,1}}}/IB_{y_{_{1,1}}})/(g,h)(B_{y_{_{1,1}}}/IB_{y_{_{1,1}}})= \dim B_{y_{_{1,1}}}/IB_{y_{_{1,1}}}-2$$
and since $B_{y_{_{1,1}}}/IB_{y_{_{1,1}}}$ is Cohen--Macaulay, this proves the claim.

\smallskip

Summing up  we have shown that $P$ has codimension at least $m(m-2)-2+2=m(m-2)$.
Therefore, $\dim V(f)^*=m^2-2-{\rm cod} \,V(f)^* \leq m^2-2-m(m-2)=2m-2$, as was to be shown.
\qed

\begin{Theorem}\label{expected_mult}
	Let $f=\det(\mathcal{GC})$. Then $f$ is a factor of its Hessian determinant $h(f)$ and has the expected multiplicity as such.
\end{Theorem}
\demo
Fix the notation of Theorem~\ref{primality_generic_cloned} (iv)  and its proof in the next subsection.
As seen there, the identity
$${\rm adj}({\rm adj}(\mathcal{GC}))=\det(\mathcal{GC})^{m-2}\cdot\mathcal{GC}$$
shows that the inverse map to the birational map defined by the minors $\boldsymbol{\Delta}$ has $\widetilde{\mathbb{D}}:=\{\mathbb{D}_{i,j}\,|\, (i,j)\neq (m,m)\}$ as its set of coordinates.
By (\ref{composite_is_identity}) one has the equality
\begin{equation}\label{composition}
\widetilde{\mathbb{D}}\circ\boldsymbol\Delta(\xx)=f^{m-2}\cdot(\xx).
\end{equation}
Applying the chain rule to (\ref{composition}) yields
\begin{equation}\label{chainrule}
\Theta(\widetilde{\mathbb{D}})(\boldsymbol\Delta)\cdot \Theta(\boldsymbol\Delta)=f^{m-2}\cdot\mathcal{I}+\frac{1}{m-2}f^{m-3}\cdot(\xx)^t\cdot{\rm Grad}(f),
\end{equation}
where $\Theta(S)$ denotes the Jacobian matrix of a set $S$ of polynomials,  ${\rm Grad}(f)$ stands for the row vector of the partial derivatives of $f$ and  $\mathcal{I}$ is the identity matrix of order $m-1$.

Write
$$\Theta(\widetilde{\mathbb{D}})(\boldsymbol\Delta)=
\left(\begin{array}{c|c}
A&U
\end{array}\right),
\quad
\Theta(\boldsymbol\Delta)=
\left(\begin{array}{c}
B\\
\hline
V
\end{array}\right),
$$
where $U$ designates the column vector of the partial derivatives of the elements of  $\widetilde{\mathbb{D}}$ with respect to $y_{m,m}$  further evaluated at $\boldsymbol\Delta$, while $V$ stands for the row vector of the partial derivatives of $\Delta_{m,m}$ with respect to the $\xx$-variables. 

Letting $E$ denote the elementary matrix obtained from $\mathcal I$ by adding the $m$th row to the $(m-1)$th row, one further has
\begin{equation*}
\Theta(\widetilde{\mathbb{D}})(\boldsymbol\Delta)\cdot E^{-1}=
\left(\begin{array}{c|c}
A'&U
\end{array}\right),
\quad
E\cdot \Theta(\boldsymbol\Delta)=
\left(\begin{array}{c}
H_f\\
\hline \\ [-12pt]
V
\end{array}\right)
\end{equation*}
where $H_f$ stands for the Hessian matrix of $f$  with respect to the $\xx$-variables.

Applying these values to (\ref{chainrule}) obtains
\begin{equation}\label{alinhahf}
A'\cdot H_f= f^{m-2}\cdot\mathcal{I}+\mathcal{A}
\end{equation}
where $\mathcal{A}=\frac{1}{m-2}f^{m-3}\cdot(\xx)^t\cdot{\rm Grad}(f)-U\cdot V$.

On the other hand, note that $\left(\partial{\mathbb{D}_{ij}}/\partial{y_{r,s}}\right)(\boldsymbol\Delta)$ 
is the $(m-2)$-minor of adj$(\mathcal{GC})$ ommiting the $i$th and $r$th rows and the  $j$th and $s$th columns.
Therefore, by a classical identity 

\begin{equation}\label{identidadeJacob}
\left(\partial{\mathbb{D}_{ij}}/\partial{y_{r,s}}\right)(\boldsymbol\Delta)= f^{m-3}(\Delta_{i,j}\Delta_{r,s}-\Delta_{i,s}\Delta_{r,j}).
\end{equation}
As a consequence,  $U=f^{m-3}\cdot U'$ where the entries of $U'$ are certain $2$-minors of  $\mathcal{GC}$. 
We now get

\begin{equation*}
A'\cdot H_f=f^{m-3}( f\cdot\mathcal{I}+\mathcal{A}')
\end{equation*}
where $\mathcal{A}'=\frac{1}{m-2}\cdot(\xx)^t\cdot{\rm Grad}(f)-U'\cdot V.$
Thus,

\begin{equation}\label{detproduto}
\det(A')\cdot\det(H_f)=f^{(m-3)(m^2-1)}\det(f\cdot\mathcal{I}+\mathcal{A}')
\end{equation}

Set $n:=m^2-1,$ $[n]:=\{1,\ldots,n\}$ and let $\Delta_{[n]\setminus\{i_1,\ldots,i_{k}\}}$ denote the principal $(n-k)$-minor of $\mathcal{A}'$ with rows and columns $[n]\setminus\{i_1,\ldots,i_{k}\}$.
Also note that $\mathcal{A}'$ has rank at most $2$ since it is a sum of matrices of rank $1$.
Using a classical formula for the determinant of a sum where one of the summands is a diagonal matrix (see, e.g., \cite[Lemma 2.3]{SiVi}), one has
{\small
	\begin{eqnarray*}\nonumber
		\det\left(f\cdot\mathcal{I}+\mathcal{A}'\right)\kern-8pt&=& \kern-8pt	\det(\mathcal{A}')+f\left(\sum_{i}\Delta_{[n]-\{i\}}\right)+\ldots+f^{n-1}\left(\sum_{1\leq i_1<\ldots<i_{n-1}\leq n}\Delta_{[n]\setminus\{i_1,\ldots,i_{n-1}\}}\right)+f^{n}\\ \nonumber
		\kern-8pt&=& \kern-8pt f^{n-2}\left(\sum_{1\leq i_1<\ldots<i_{n-2}\leq n}\Delta_{[n]\setminus\{i_1,\ldots,i_{n-2}\}}\right)+f^{n-1}\cdot{\rm \mbox{trace}}(\mathcal{A}')+f^{n}.
\end{eqnarray*}}

Setting $G:=\displaystyle\sum_{1\leq i_1<\ldots<i_{n-2}\leq n}\Delta_{[n]\setminus\{i_1,\ldots,i_{n-2}\}}+f\cdot{\rm \mbox{trace}}(\mathcal{A}')+f^{2}$ and substituting in (\ref{detproduto}) gives
\begin{equation}\label{potgrande}
\det(A')\cdot\det(H_f)=f^{(m-3)(m^2-1)+m^2-3}\cdot G=f^{(m^2-1)(m-2)-2}\cdot G
\end{equation}

Suppose for a moment that $G\neq 0$. In this case, clearly $\det(A')\neq 0$ and a degree argument shows that some positive power of $f$ divides $h(f)$.
Indeed, by construction, 
$$\deg(\det(A'))=(m^2-1)(m(m-3)+2)<m((m^2-1)(m-2)-2)= \deg(f^{(m^2-1)(m-2)-2}).$$
Therefore, we are left with proving that $G$ does not vanish.

Note that the vanishing of $G$ would imply in particular that $f$ is an integral element over the $k$-subalgebra generated by the entries of $\mathcal{A}'$. This would possibly be forbidden if the latter could be proved to be integrally closed.
Due to the difficulty of this verification we resort to a direct inspection.
For this, note that 
\begin{eqnarray}
\mbox{trace}(\mathcal{A}')&=&\frac{1}{m-2}\,\mbox{trace}((\xx)^t\cdot{\rm Grad}(f))-\mbox{trace}(U'V)\nonumber\\
&=&\frac{m}{m-2}\, f-\sum_{ij\atop i,j\neq m} (a_{ij}a_{m,m}-a_{i,m}a_{mj})\frac{\partial\Delta_{m,m}}{\partial x_{ij}}\nonumber\\
&=&\frac{m}{m-2}\, f-\left(a_{mm}\sum_{ij} a_{ij}\frac{\partial\Delta_{m,m}}{\partial x_{ij}}-\sum_{ij\atop i,j\neq m}a_{im}a_{mj}\frac{\partial\Delta_{m,m}}{\partial x_{ij}}\right)\nonumber\\
&=&\frac{m}{m-2}\, f-(m-1)a_{m,m}\Delta_{m,m}+\sum_{ij\atop i,j\neq m}a_{im}a_{mj}\frac{\partial\Delta_{m,m}}{\partial x_{ij}}.\nonumber
\end{eqnarray}
Setting  $q_{ij}:=a_{ij}a_{mm}-a_{im}a_{mj}$, one has

{\small \begin{eqnarray}
	\displaystyle\sum_{1\leq i_1<\ldots<i_{n-2}\leq n}\Delta_{[n]\setminus\{i_1,\ldots,i_{n-2}\}}
	&=&
	\sum_{(ij),(rs)}\det\left(\begin{array}{cc}
	x_{ij}\frac{\partial f}{\partial x_{ij}}-q_{ij}\frac{\partial\Delta_{mm}}{\partial x_{ij}}& x_{ij}\frac{\partial f}{\partial x_{rs}}-q_{ij}\frac{\partial\Delta_{mm}}{\partial x_{rs}}\nonumber\\[5pt]
	x_{rs}\frac{\partial f}{\partial x_{ij}}-q_{rs}\frac{\partial\Delta_{mm}}{\partial x_{ij}}& x_{rs}\frac{\partial f}{\partial x_{rs}}-q_{rs}\frac{\partial\Delta_{mm}}{\partial x_{rs}}
	\end{array}\right)\nonumber\\
	&=&
	\sum_{(ij),(rs)}\det\left(\begin{array}{cc}
	x_{ij}&-q_{ij}\\
	x_{rs}&-q_{rs}\end{array}\right)\cdot \det\left(\begin{array}{cc}
	\frac{\partial f}{\partial x_{ij}}&\frac{\partial f}{\partial x_{rs}}\\ [5pt]
	\frac{\partial \Delta_{mm}}{\partial x_{ij}}&\frac{\partial \Delta_{mm}}{\partial x_{rs}}\end{array}\right)\nonumber\\
	&=&\sum_{(ij),(rs)}\det\left(\begin{array}{cc}
	x_{ij}&a_{im}a_{mj}\\
	x_{rs}&a_{rm}a_{ms}\end{array}\right)\cdot \det\left(\begin{array}{cc}
	\frac{\partial f}{\partial x_{ij}}&\frac{\partial f}{\partial x_{rs}}\nonumber\\[5pt]
	\frac{\partial \Delta_{mm}}{\partial x_{ij}}&\frac{\partial \Delta_{mm}}{\partial x_{rs}}\end{array}\right)\nonumber
	\end{eqnarray}}

Thus, we can write $G=G_1+G_2$, where 
$$G_1=f\left(\sum_{ij\atop i,j\neq 0} a_{im}a_{mj}\frac{\partial\Delta_{m,m}}{\partial x_{ij}}\right)+\sum_{(ij),(rs)}\det\left(\begin{array}{cc}
x_{ij}&a_{im}a_{mj}\\
x_{rs}&a_{rm}a_{ms}\end{array}\right)\cdot \det\left(\begin{array}{cc}
\frac{\partial f}{\partial x_{ij}}&\frac{\partial f}{\partial x_{rs}}\\
\frac{\partial \Delta_{mm}}{\partial x_{ij}}&\frac{\partial \Delta_{mm}}{\partial x_{rs}}\end{array}\right)$$
and 
$$G_2= f\left(\frac{m}{m-2}\cdot f-(m-1)a_{m,m}\Delta_{m,m}\right)+f^2.$$
Inspecting the summands of $G_1$ and $G_2$ one sees that the degree of $x_{m-1,m-1}$ in $G_1$ is at most $3$, while that of $x_{m-1,m-1}$ in $G_2$ is $4$.
This shows that $G\neq 0$.
\qed

\begin{Remark}\rm
	(1) The factor $h(f)/ f^{m(m-2)-1}$ coincides with the determinant of the $2\times 2$ submatrix with rows $m-1,m$ and columns $m-1,m$.
	
	(2) There is a mistaken assertion in the proof of \cite[Proposition 3.2 (a)]{MAron} to the effect that the dual variety to the generic determinant $f$ is the variety of submaximal minors. This is of course nonsense since the dual variety is the variety of the $2\times 2$ minors. This wrong assertion in loc. cit. actually serves no purpose in the proof since (a) follows simply from the cofactor relation. At the other end, this nonsense reflects on the second assertion of item (c) of the same proposition which is therefore also flawed. The right conclusion is that the multiplicity of $f$ as a factor of its Hessian has indeed the expected multiplicity, since the variety of the $2\times 2$ minors has codimension $(m-1)^2$ and hence $m^2-1-(m-1)^2=m(m-2)$ as desired.
	Likewise, \cite[Conjecture 3.4 (a)]{MAron} should be read as affirmative without exception.
	
	(3) We have been kindly informed by J. Landsberg that either Theorem~\ref{dim_dual} or Theorem~\ref{expected_mult}, or perhaps both -- which were stated as a conjecture in the first version of this prepint posted on the arXiv --  have been obtained in \cite{Landsberg} by geometric means.
	For our misfortune, we were not able to trace in the mentioned work the precise statements expressing the above contents.
\end{Remark}

\section{Degeneration by zeros}\label{zeros}

Recall from the previous section the cloning degeneration where an entry is cloned along the same row or column of the original generic matrix.
As mentioned before, up to elementary operations of rows and/or columns the resulting matrix has a zero entry.
A glimpse of this first status has been tackled in \cite[Proposition 4.9 (a)]{MAron}.

This procedure can be repeated to add more zeros.
Aiming at a uniform treatment of all these cases, we will fix integers $m,r$ with $1\leq r\leq m-2$ and consider the following degeneration of the $m\times m$ generic matrix:

{\scriptsize
	\begin{equation}\label{generic-zeros}
	\left(
	\begin{array}{cccccccc}
	x_{1,1}&\ldots & x_{1,m-r}& x_{1,m-r+1} &  x_{1,m-r+2}&\ldots & x_{1,m-1} & x_{1,m}\\
	\vdots & \ldots & \vdots & \vdots & \vdots & \ldots & \vdots &\vdots \\
	x_{m-r,1}& \ldots & x_{m-r,m-r}& x_{m-r,m-r+1} & x_{m-r,m-r+2}&\ldots & x_{m-r,m-1} & x_{m-r,m}\\
	x_{m-r+1,1}&\ldots & x_{m-r+1,m-r}& x_{m-r+1,m-r+1} & x_{m-r+1,m-r+2} &\ldots & x_{m-r+1,m-1} & 0\\
	x_{m-r+2,1}&\ldots & x_{m-r+2,m-r}& x_{m-r+2,m-r+1} & x_{m-r+2,m-r+2} &\ldots & 0 & 0\\
	\vdots & \ldots & \vdots & \vdots & \vdots&\iddots &\vdots & \vdots\\
	x_{m-1,1} &\ldots & x_{m-1,m-r}&  x_{m-1,m-r+1}& 0&\ldots & 0 & 0\\[4pt]
	x_{m,1} &\ldots & x_{m,m-r}& 0& 0&\ldots & 0 & 0\\ [3pt]
	\end{array}
	\right)
	\end{equation}
}

Assuming $m$ is fixed in the context, let us denote the above matrix by $\mathcal{DG}(r)$. 

\subsection{Polar behavior}

\begin{Theorem}\label{polar-zeros}
Let $R=k[\xx]$ denote the polynomial ring in the nonzero entries of $\mathcal{DG}(r)$, with $1\leq r\leq m-2$, let $f :=\det \mathcal{DG}(r)$ and let $J\subset R$ denote the gradient ideal of $f$. Then:
\begin{itemize}
	\item[{\rm (a)}] $f$ is irreducible.
	\item[{\rm (b)}] $J$ has maximal linear rank.
	\item[{\rm (c)}] The homogeneous coordinate ring  of the polar variety of $f$ in $\pp^{m^2-{{r+1}\choose {2}}-1}$ is a Gorenstein ladder determinantal ring of dimension $m^2-r(r+1);$ in particular, the analytic spread of $J$ is $m^2-r(r+1)$.
\end{itemize}
\end{Theorem}
\demo
(a) 
Expanding the determinant by Laplace along the first row, we can write $f=x_{1,1}\Delta_{1,1}+g$, where $\Delta_{1,1}$ is the cofactor of $x_{1,1}$.
Clearly, both $\Delta_{1,1}$ and $g$ belong to the polynomial subring omitting the variable $x_{1,1}$.
Thus, in order to show that $f$ is irreducible it suffices to prove that it is a primitive polynomial (of degree $1$) in $k[x_{1,2},\ldots, x_{m,m-r}][x_{1,1}]$.
In other words, we need to check that no irreducible factor of $\Delta_{1,1}$ is a factor of $g$.

We induct on $m\geq r+2$. 
If $m=r+2$ then $\Delta_{1,1}=x_{2,m} x_{3,m-1} \cdots x_{m-1,3} x_{m,2}$, while the initial term of $g$ in the revlex monomial order is
$$
	{\rm in}(g) =	{\rm in}(f)=x_{1,m} x_{2,m-1} \cdots  x_{m,1}.
$$
Thus, assume that $m> r+2$. 
By the inductive step, $\Delta_{1,1}$ is irreducible being the determinant of an $(m-1)\times (m-1)$ matrix of the same kind (same $r$).
But $\deg(\Delta_{1,1})= \deg(g)-1$.
Therefore, it suffices to show that $\Delta_{1,1}$ is not a factor of $g$.
Supposing it were, we would get that $f$ is multiple of $\Delta_{1,1}$ by a linear factor -- this is clearly impossible.

Once more, an alternative  argument is to use that the ideal $J$  has codimension $4$, as will be shown independently in Theorem~\ref{primality_generic_zeros} (b). 
Therefore, the ring $R/(f)$ is locally regular in codimension at least one, so it must be normal.
But $f$ is homogeneous, hence irreducible.

\medskip

(b)  
The proof is similar to the one of Theorem~\ref{cloning_generic} (iii), but there is a numerical diversion and, besides, the cases where $r> m-r-1$ and $r\leq m-r-1$ keep slight differences.

Let $f_{i,j}$ denote the  $x_{i,j}$-derivative of $f$ and let $\Delta_{j,i}$ stand for the (signed) cofactor of $x_{i,j}$ on $\mathcal{DG}(r)$.
We first assume that  $r> m-r-1$. The Cauchy cofactor formula
$$\mathcal{DG}(r) \cdot {\rm adj}(\mathcal{DG}(r))={\rm adj}(\mathcal{DG}(r))\cdot \mathcal{DG}(r)=\det(\mathcal{DG}(r))\,\mathbb{I}_m$$
yields by expansion the following three blocks of linear  relations involving the  (signed) cofactors of  $\mathcal{DG}(r)$:

\begin{equation}
\label{eqs1}
\left\{\begin{array}{ll}
\sum_{j=1}^m x_{i,j}\Delta_{j,k}=0 & \kern-10pt\mbox{for $1\leq i \leq m-r,\; 1\leq k  \leq m-r\;( k\neq i)$}\\ [2pt]
\sum_{j=1}^{m-l}x_{m-r+l,j}\Delta_{j,k}=0 & \mbox{for $1\leq l \leq r, \; 1\leq k  \leq m-r$}\\ [2pt] 
\sum_{j=1}^m x_{i,j}\Delta_{j,i}-\sum_{j=1}^m x_{i+1,j}\Delta_{j,i+1}=0 & \mbox{for $1\leq i\leq m-r-1$}
\end{array}
\right.
\end{equation}
with $m^2-rm-1$ such relations;

\begin{equation}
\label{eqs2}\left\{\begin{array}{ll}
\sum_{i=1}^{m} x_{i,j}\Delta_{k,i}=0 & \mbox{for $1\leq j \leq m-r, \; 1\leq k  \leq m-r\;( k\neq j)$}\\ [3pt]
\sum_{i=1}^{m-l} x_{i,m-r+l}\Delta_{k,i}=0 & \mbox{for $1\leq l \leq 2r-m+1, \;  1\leq k  \leq m-r$}
\end{array}
\right.
\end{equation}
with    $(m-r)(m-r-1)+(m-r)(2r-m+1)=r(m-r)$ such relations; and

\begin{equation}
\label{eqs3}
\left\{\begin{array}{ll}
\sum_{i=1}^{m-l} x_{i,m-r+l}\Delta_{m-r+l,i}-\sum_{j=1}^m x_{i,1}\Delta_{1,i}=0 & \mbox{for $1\leq l\leq r-1$}\\ [3pt]
\sum_{i=1}^{m-k} x_{i,m-r+k}\Delta_{m-r+l,i}=0 & \kern-10pt\mbox{for $1\leq l \leq r-2, \;l+1\leq k\leq r-1$}
\end{array}
\right.
\end{equation}
with $r(r-1)/2$ such relations. 

\smallskip

Similarly, when $r\leq m-r-1$, the classical Cauchy cofactor formula outputs by expansion  three blocks of linear relations involving the (signed) cofactors of $\mathcal{DG}$. Here, the first and third blocks are, respectively,  exactly as the above ones, while the second one requires a modification due to the inequality reversal; namely, we get

\begin{equation}\label{eqs'}
\left\{\begin{array}{ll}
\sum_{i=1}^{m} x_{i,j}\Delta_{k,i}=0 & \mbox{for $1\leq j \leq r+1, \;1\leq k  \leq r\;( k\neq j)$}\\ [3pt]
\sum_{i=1}^{m} x_{i,j}\Delta_{k,i}=0 & \mbox{for $1\leq j\leq r, \; r+1\leq k  \leq m-r$},
\end{array}
\right.
\end{equation}
with $r(m-r)$ such relations (as before).

\smallskip

Since $f_{i,j}$ coincides with the signed cofactor $\Delta_{j,i}$, any of the above relations gives a linear syzygy of the partial derivatives of $f$. Thus one has a total of $m^2-rm-1+r(m-r)+r(r-1)/2=m^2-{{r+1}\choose {2}}-1$ linear syzygies of $J$.

It remains to show that these are independent.

For this, we adopt the same strategy as in the proof of Theorem~\ref{cloning_generic} (iii), whereby we list the partial derivatives according to the following ordering of the nonzero entries:  we traverse the first row from left to right, then the second row in the same way, and so on until we reach the last row with no zero entry; thereafter we start from the first row having a zero and travel along the columns, from left to right, on each column from top to bottom, till we all nonzero entries are counted.

Thus, the desired ordering is depicted in the following scheme, where we once more used arrows for easy reading:
\begin{eqnarray*} \nonumber
	\lefteqn{ \kern-40pt x_{1,1},x_{1,2},\ldots, x_{1,m}\rightsquigarrow x_{2,1}, x_{2,2},\ldots, x_{2,m} \rightsquigarrow \ldots \rightsquigarrow x_{m-r,1},x_{m-r,2}\ldots, x_{m-r,m} \rightsquigarrow }  \kern-40pt\\ 
	&&\kern-20pt x_{m-r+1,1},   \ldots, x_{m,1} \rightsquigarrow x_{m-r+1,2},\ldots ,x_{m,2}
	\rightsquigarrow \ldots \rightsquigarrow x_{m-r+1,m-r},x_{m-r+2,m-r},\ldots, x_{m,m-r} \\
	&& \kern-20pt\rightsquigarrow x_{m-r+1,m-r+1}, \ldots, x_{m-1,m-r+1} \rightsquigarrow x_{m-r+1,m-r+2}, \ldots,x_{m-2,m-r+2 }\\ 
	&&\kern-20pt \rightsquigarrow \ldots \rightsquigarrow x_{m-r+1,m-2}, x_{m-r+2,m-2}\rightsquigarrow x_{m-r+1,m-1} 
\end{eqnarray*}
With this ordering the above linear relations translate into linear syzygies collected in the following block matrix

{\small
	$$\mathcal{M}=\left(\begin{array}{cccc|cccc|ccccc}
	\varphi_1 &   \ldots &         &              &                 &     &        &                                     &            &     &      &\\
	\mathbf{0}       &\varphi_2 &\ldots   &              &                 &     &        &                                     &            &    &    &\\
	\vdots    &\vdots    & \ddots  & \vdots       &                 &     &        &                                   &            &     &   &\\
	\mathbf{0}        & \mathbf{0}        & \ldots  & \varphi_{m-r}&                 &    &        &                                    &            &       & &\\ 
	\hline\\ [-10pt]
	\mathbf{0}_r^{m-1} & \mathbf{0}_r^m & \ldots  &  \mathbf{0}_r^m	      & \varphi_{r}^1   &     &        &                       &                     &        & &\\
	\mathbf{0}_r^{m-1} & \mathbf{0}_r^m  &  \ldots & \mathbf{0}_r^m        &    \mathbf{0}_r^r       & \varphi_{r}^2   &  &                       &                  &        & &\\
	\vdots    & \vdots & \ldots  & \vdots       &  \vdots          &  \vdots        &\ddots & & & & &\\
	\mathbf{0}_r^{m-1} & \mathbf{0}_r^m& \ldots  &  \mathbf{0}_r^m	      &    \mathbf{0}_r^r       &  \mathbf{0}_r^r   & \ldots &   \varphi_{r}^{(m-r)}&                 &        & &\\
	\hline\\ [-10pt]
	\mathbf{0}_{r-1}^{m-1}    & \mathbf{0}_{r-1}^{m}        & \ldots  &     \mathbf{0}_{r-1}^m  &    \mathbf{0}_{r-1}^r          & \mathbf{0}_{r-1}^r      &    \ldots   &  \mathbf{0}_{r-1}^r    & \Phi_1 &  &  &  \\
	\mathbf{0}_ {r-2}^{m-1}   & \mathbf{0}_{r-2}^m        & \ldots  &     \mathbf{0}_{r-2}^m    &    \mathbf{0}_{r-2}^r                & \mathbf{0}_{r-2}^r      &          \ldots   &  \mathbf{0}_{r-2}^r  &  \mathbf{0}_{r-2}^{r-1}  & \Phi_2    & &  \\
	\vdots   & \vdots        & \vdots  &     \vdots    &    \vdots                & \vdots      &          \vdots  &  \vdots &  \vdots  & \vdots   &\ddots &  \\
	\mathbf{0}_1^{m-1}    & \mathbf{0}_1^m       & \ldots  & \mathbf{0}_1^m         & \mathbf{0}_1^r                   & \mathbf{0}_1^r     &         \ldots    &  \mathbf{0}_1^r    & \mathbf{0}_1^{r-1}   & \mathbf{0}_1^{r-2}  & \ldots &\Phi_{r-1} \\
	\end{array}
	\right),$$
}
where:

\begin{itemize}
	\item $\varphi_1$ is the matrix obtained  from the  transpose $\mathcal{DG}(r)^t$ of $\mathcal{DG}(r)$ by omitting the first column
	\item $\varphi_2,\ldots,\varphi_{m-r}$ are each a copy of  $\mathcal{DG}(r)^t$ (up to column permutation);
	\item    When $r>m-r-1$, setting $\mathfrak{d}=2r-m+1$, for $i=1,\ldots,m-r$ one has that $\varphi_{r}^{i}$ is the $r\times r$ minor omitting the $i$th column  of the following submatrix  of $\mathcal{DG}(r)$:

	$$\kern-7pt\left( \begin{matrix}
	x_{m-r+1,1}& \ldots & x_{m-r+1,m-r}& x_{m-r+1,m-r+1}& x_{m-r+1,m-r+2} 	&\ldots& x_{m-r+1,r+1}\\ 
	\vdots & \vdots &\vdots & \vdots & \vdots 	& \vdots & \vdots \\
	x_{m-\mathfrak{d},1}& \ldots & x_{m-\mathfrak{d},m-r}& x_{m-\mathfrak{d},m-r+1}& x_{m-\mathfrak{d},m-r+2} &\ldots& x_{m-\mathfrak{d},r+1}\\
	x_{m-\mathfrak{d}+1,1}& \ldots & x_{m-\mathfrak{d}+1,m-r}& x_{m-\mathfrak{d}+1,m-r+1}& x_{m-\mathfrak{d}+1,m-r+2} &\ldots& 0\\
	\vdots & \vdots &\vdots & \vdots & \vdots 	& \vdots & \vdots \\ 
	x_{m-1,1}& \ldots & x_{m-1,m-r}& x_{m-1,m-r+1}& 0 &\ldots& 0\\
	x_{m,1}& \ldots & x_{m,m-r}& 0& 0 &\ldots& 0\\
	\end{matrix}\right). $$
	
	When $r\leq m-r-1$, consider the  following submatrix of $\mathcal{DG}(r)$:
	{\Large	
		$$\left( \begin{matrix}
		x_{m-r+1,1}& \ldots & x_{m-r+1,r}& x_{m-r+1,r+1}\\
		\vdots & \vdots &\vdots &  \vdots \\
		_{m,1}& \ldots & x_{m, r} & x_{m,r+1}		\end{matrix}\right). $$
	}
Then, for  $i=1,\ldots,r$ (respectively,  for $i=r+1,\ldots, m-r$) $\varphi_{r}^{i}$ denotes the $r\times r$ submatrix obtained  by omitting the $i$th column (respectively, the last column).

	\item Each $\mathbf{0}$ under $\varphi_1$  is an $m\times (m-1)$  block of zeros and each $\mathbf{0}$ under $\varphi_i$ is an $m\times m$  block of zeros for $i=2,\ldots,m-r-1$ ;
	\item $\mathbf{0}_l^c$ denotes  an $l\times c $ block of zeros.
	\item $\Phi_i$ is the following $(r-i)\times (r-i)$ submatrix of $\mathcal{DG}(r)$: 
	{\Large
	$$\Phi_i=\left(\begin{matrix}
	x_{m-r+1,m-r+i} & x_{m-r+1,m-r+i+1} & \ldots & x_{m-r+1,m-1}\\
	x_{m-r+2,m-r+i} & x_{m-r+2,m-r+i+1} & \ldots & 0\\
	\vdots & \vdots & \vdots & \vdots \\
	x_{m-i,m-r+i} & 0 & \ldots & 0 	
	\end{matrix}\right). $$
}
	
\end{itemize}

Next we justify why these blocks make up (linear) syzygies.  As already explained, the relations in (\ref{eqs1}), (\ref{eqs2}) and (\ref{eqs3}) yield linear syzygies of the partial derivatives of $f$. Setting $k=1$ in the first two relations of (\ref{eqs1}), the latter can be written  as $\sum_{j=1}^m x_{i,j}f_{1,j}=0$,  for $i=2,\ldots,  m-r$, and  $\sum_{j=1}^{m-l}x_{m-r+l,j}f_{1,j}=0$, for all  $l =1,\ldots, r$. By ordering the set of partial derivatives $f_{i,j}$ as explained before, the coefficients of these relations become the entries of the submatrix $\phi_1$ of $\mathcal{DG}(r)^t$ obtained by omitting its first column, as mentioned above, namely: 

{\small 
	$$\left(
	\begin{array}{cccccccc}
	x_{2,1}&\ldots & x_{m-r,1}& x_{m-r+1,1} &  x_{m-r+2,1}&\ldots & x_{m-1,1} & x_{m,1}\\
	\vdots & \ldots & \vdots & \vdots & \vdots & \ldots & \vdots &\vdots \\
	x_{2,m-r}& \ldots & x_{m-r,m-r}& x_{m-r+1,m-r} & x_{m-r+2, m-r}&\ldots & x_{m-1, m-r} & x_{m,m-r}\\
	x_{2,m-r+1}&\ldots & x_{m-r,m-r+1}& x_{m-r+1,m-r+1} & x_{m-r+2, m-r+1} &\ldots & x_{m-1, m-r+1} & 0\\
	x_{2, m-r+2}&\ldots & x_{m-r,m-r+2}& x_{m-r+1,m-r+2} & x_{m-r+2,m-r+2} &\ldots & 0 & 0\\
	\vdots & \ldots & \vdots & \vdots & \vdots&\iddots &\vdots & \vdots\\
	x_{2,m-1} &\ldots & x_{m-r,m-1}&  x_{m-r+1,m-1}& 0&\ldots & 0 & 0\\
	x_{2,m} &\ldots & x_{m-r,m}& 0& 0&\ldots & 0 & 0
	\end{array}
	\right)$$}

Getting $\varphi_k$, for $k=2,\ldots , m-r$, is similar, namely, we use  again the first two relations  in the block (\ref{eqs1}) retrieving the submatrix of $\mathcal{DG}(r)^t$ excluding the $k$th column and replacing  it with an extra  column that comes  from the last relation in (\ref{eqs1}) by taking $i=k-1$.  

Continuing, for each $i=1,\ldots, m-r$ the block  $\varphi_r^i$  comes from the relations in the blocks (\ref{eqs2}), if  $r>m-r-1$, or  (\ref{eqs'}),  if $r\leq m-r-1$, by setting $k=i$. Finally, for each $i=1,\ldots, r-1$, the block $\Phi_i$ comes from the relations in (\ref{eqs3}) by setting $l=i$.

This proves the claim about the large matrix above.
Counting through the sizes of the various blocks, one sees that this matrix is $(m^2-{{r+1}\choose {2}})\times (m^2-{{r+1}\choose {2}}-1)$.
Omitting its first row obtains a square block-diagonal submatrix  where each block has nonzero determinant. Thus, the linear rank of $J$ attains the maximum.

\smallskip

(c) Note that the polar map can be thought as the map of $\pp^{m^2-{{r+1}\choose {2}}-1}$ to itself defined by the partial derivatives of $f$.
As such, the polar variety will be described in terms of defining equations in the original $\xx$-variables.

\begin{figure}[h]
	\centering
	\includegraphics[width=14cm, height=7cm]{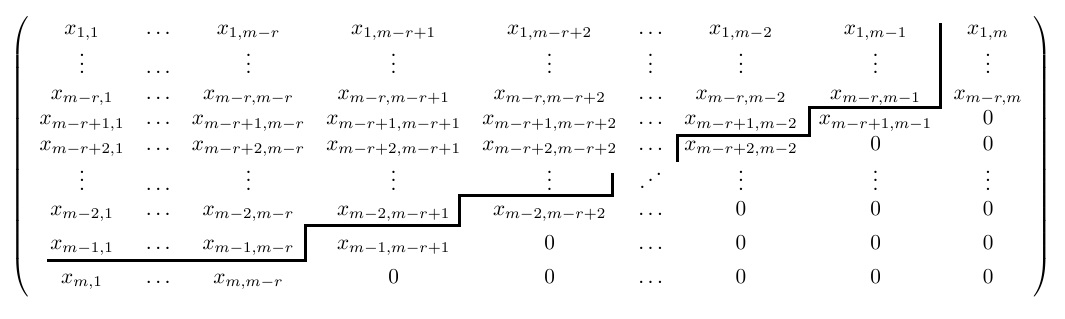}\\
	\caption{stair-like polygonal.}\label{r-stair-like}
\end{figure}

Let $\mathcal{L}=\mathcal{L}(m,r)$ denote  the set of variables in $\mathcal{DG}(r)$ lying to the left and above the stair-like polygonal in  Figure~\ref{r-stair-like}  and let $I_{m-r}(\mathcal{L})$ stand for the ideal generated by the $(m-r)\times (m-r)$ minors of $\mathcal{DG}(r)$ involving only the variables in $\mathcal{L}$.  

Since $\mathcal{L}$ can be extended to a fully generic matrix of size $(m-1)\times (m-1)$, the ring $K[\mathcal{L}]/I_{m-r}(\mathcal{L}) $ is one of the so-called ladder determinantal  rings. 

{\sc Claim:} The homogeneous defining ideal of the image of the polar map of $f$ contains the ideal $I_{m-r}(\mathcal{L})$. 

Let $x_{i,j}$ denote a nonzero entry of $\mathcal{DG}(r)$. Since the nonzero entries of the matrix are independent variables, it follows easily from the Laplace expansion along the $i$th row that the $x_{i,j}$-derivative $f_{i,j}$ of $f$ coincides with the (signed) cofactor of  $x_{i,j}$, heretofore denoted $\Delta_{j,i}$. 

Given  integers $1\leq i_1< i_2 <\ldots <i_{m-r}\leq m-1$, consider the following submatrix of the transpose matrix of cofactors:

$$F= \left( \begin{matrix}
\Delta_{i_1,1} &\Delta_{i_1,2}& \Delta_{i_1,3}&\cdots & \Delta_{i_1,m-i_{m-r}+(m-r-1)}\\
\Delta_{i_2,1} &\Delta_{i_2,2}& \Delta_{i_2,3}&\cdots & \Delta_{i_2,m-i_{m-r}+(m-r-1)}\\
\vdots  &\vdots & \vdots &\cdots & \vdots \\
\Delta_{i_{m-r},1} &\Delta_{i_{m-r},2}& \Delta_{i_{m-r},3}&\cdots & \Delta_{i_{m-r},m-i_{m-r}+(m-r-1)}
\end{matrix}\right).$$

Letting 
$$  C=\left(  \begin{matrix}
x_{1,i_{m-r}+1} & x_{1,i_{m-r}+2}& \cdots & x_{1,m-1} & x_{1,m}\\
\vdots & \vdots & \cdots & \vdots & \vdots \\
x_{m-r,i_{m-r}+1} & x_{m-r,i_{m-r}+2}& \cdots & x_{m-r,m-1} & x_{m-r,m}\\
x_{m-r+1,i_{m-r}+1} & x_{m-r+1,i_{m-r}+2}& \cdots & x_{m-r+1,m-1} &    0\\
\vdots & \vdots & \cdots &\vdots  &  \vdots\\
x_{m-i_{m-r}+(m-r-2),i_{m-r}+1} & x_{m-i_{m-r}+(m-r-2),i_{m-r}+2}& \cdots & 0 &    0\\
x_{m-i_{m-r}+(m-r-1),i_{m-r}+1} & 0 & \cdots & 0 &    0\\
\end{matrix}\right),$$
the cofactor identity ${\rm adj}( \mathcal{DG}(r))\cdot \mathcal{DG}(r) = \det(\mathcal{DG}(r))\mathbb{I}_m$ yields the relation
$$F\cdot C=0.$$

Since the  columns of $C$ are linearly independent, it follows that the rank of $F$ is at most $m-i_{m-r}+(m-r-1) -(m-i_{m-r})=(m-r)-1$. In other words,
the maximal minors of the following matrix 
$$\left( \begin{matrix}
x_{i_1,1} &x_{i_1,2}& x_{i_1,3}&\cdots & x_{i_1,m-i_{m-r}+(m-r-1)}\\
x_{i_2,1} &x_{i_2,2}& x_{i_2,3}&\cdots & x_{i_2,m-i_{m-r}+(m-r-1)}\\
\vdots  &\vdots & \vdots &\cdots & \vdots \\
x_{i_{m-r},1} &x_{i_{m-r},2}& x_{i_{m-r},3}&\cdots & x_{i_{m-r},m-i_{m-r}+(m-r-1)}
\end{matrix}\right).$$
all vanish on the partial derivatives of $f$, thus proving the claim.

\medskip

{\sc Claim:} The codimension of  the ideal  $I_{m-r}(\mathcal{L}(m,r))$ is at least ${{r+1}\choose {2}}$. 

Let us note that the codimension of this ladder ideal could be obtained by the general principle described in \cite{HeTr} (see also \cite{ladder2}), as done in the proof of Theorem~\ref{dim_dual}.
However, in this structured situation we prefer to give an independent argument.

For this we  induct with the following inductive hypothesis: let $1\leq i\leq r-1$; then for any $(m-i)\times (m-i)$ matrix of the form $\mathcal{DG}(r-i)$, the ideal $I_{m-i-(r-i)}(\mathcal{L}(m-i,r-i))$ has codimension at least ${{r-i+1}\choose {2}}$.
Note that $m-i-(r-i)=m-r$, hence the size of the inner minors does not change in the inductive step.

We descend with regard to $i$; thus, the induction step starts out at $i=r-1$, hence $r-i=1$ and $m-i=m-(r-1)=m-r+1$ and since $m-r\geq 2$ by assumption, then $3\leq m-(r-1)\leq m-1$. Rewriting $n:=m-r+1$, we are in the situation of an $n\times n\, (n\geq 3)$ matrix of the form  $\mathcal{DG}(1)$.
Clearly, then the ladder ideal $I_{n-1}(\mathcal{L}(n-1,1))$ is a principal ideal generated by the $(n-1)\times (n-1)$ minor of $\mathcal{DG}(1)$ of the first $n-1$ rows and columns. Therefore, its codimension is $1$ as desired. 

\begin{figure}[h]
	\centering
	\includegraphics[width=15cm, height=6cm]{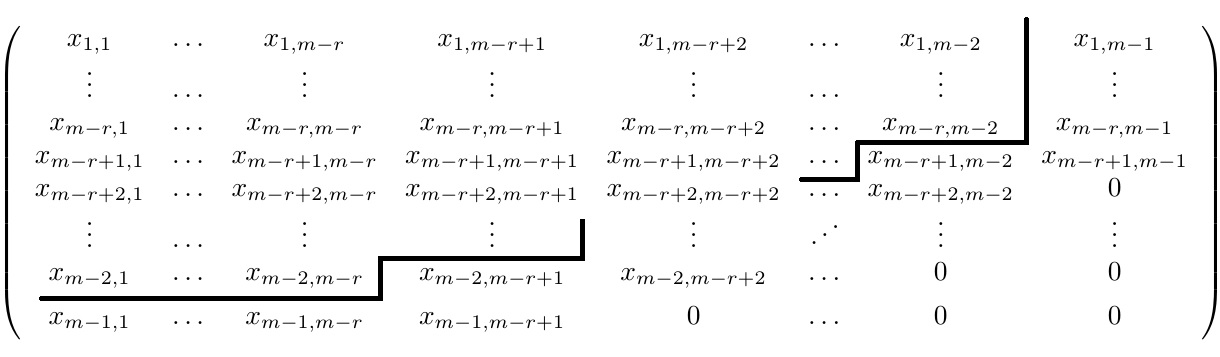}\\
	\caption{Sub-stair-like inductive.}\label{sub-stairs-inductive}
\end{figure}

To construct a suitable inductive precedent, let $\widetilde{\mathcal{L}}$ denote the set of variables  that are to the left and above  the stair-like polygonal in Figure~\ref{sub-stairs-inductive}
and denote $I_{m-r}(\mathcal{\widetilde{L}})$ the ideal generated by the $(m-r)\times (m-r)$ minors of $\mathcal{DG}(r)$ involving only the variables in $\mathcal{\widetilde{L}}$.
Note that $\mathcal{\widetilde{L}}$ is of the form $\mathcal{L}(m-1,r-1)$ relative to a matrix  of the form $\mathcal{DG}(r-1))$.
Clearly, $I_{m-r}(\mathcal{\widetilde{L}})$ it too is a ladder determinantal ideal  on a suitable $(m-2)\times (m-2)$ generic matrix; in particular, it is a Cohen-Macaulay prime ideal (see \cite{Nar} for primeness and \cite{HeTr} for Cohen--Macaulayness).  
By the inductive hypothesis, the codimension of $I_{m-r}(\mathcal{\widetilde{L}})$ is  at least ${{r}\choose {2}}$.

Note that $\widetilde{\mathcal{L}}$ is a subset of $\mathcal{L}$, hence there is a natural ring surjection
$$S:=\frac{k[\mathcal{L}]}{I_{m-r}(\mathcal{\widetilde{L}})\,k[\mathcal{L}]}=\frac{k[\mathcal{\widetilde{L}}]}{I_{m-r}(\mathcal{\widetilde{L}})}[\mathcal{L}\setminus\mathcal{\widetilde{L}}]\surjects \frac{k[\mathcal{L}]}{I_{m-r}(\mathcal{L})}
$$

Since  ${{r}\choose {2}}+r= {{r+1}\choose {2}}$, it suffices to exhibit $r$ elements of   $I_{m-r}(\mathcal{L}) $ forming a regular sequence on the ring $S:=k[\mathcal{L}]/I_{m-r}(\mathcal{\widetilde{L}})k[\mathcal{L}]$.

 Consider the matrices
\begin{equation}
\label{DELTA}
\left(
\begin{array}{ccc|c}
x_{1,1} & \ldots & x_{1, m-r-1} & x_{1,m-i}\\
\vdots & \vdots & \vdots &\vdots \\
x_{m-r-1,1}& \ldots & x_{m-r-1,m-r-1}& x_{m-r-1,m-i}\\
\hline
x_{m-r-1+i,1}& \ldots & x_{m-r-1+i,m-r-1}& x_{m-r-1+i,m-i}
\end{array}
\right)
\end{equation}
for $i=1,\ldots,r$.
Let $\Delta_i \in I_{m-r}(\mathcal{L})$ denote the determinant of the above matrix, for $i=1,\ldots,r$.

The claim is that $\Delta:=\{\Delta_1,\ldots,\Delta_{r}\}$ is a regular sequence on $S$.

Let $\delta$ denote the $(m-r-1)$-minor in  the upper left corner of (\ref{DELTA}).  Clearly, $\delta$ is a regular element on $S$ as its defining ideal is a prime ideal generated in degree $m-r$.  Therefore,  it suffices to show that the localized sequence $$\Delta_{\delta}:=\{(\Delta_1)_{\delta},\ldots,(\Delta_{r})_{\delta}\}$$  
is a regular on $S_{\delta}$. On the other hand, since $S$ is Cohen-Macaulay, it is suffices to show that ${\rm dim} \ S_{\delta}/\Delta_{\delta} S_{\delta}={\rm dim } \ S_\delta-r$.


Write $\XX':=\left\lbrace x_{m-r,m-1}, x_{m-r+1,m-2},\ldots,x_{m-2,m-r+1}, x_{m-1,m-r}\right\rbrace $.
Note that, for every $i=1,\ldots,r$, one has  $(\Delta_i)_{\delta}= x_{m-r-1+i,m-i}+ (1/\delta)\Gamma_i,$ with $x_{m-r-1+i,m-i}\in \XX' $ and $\Gamma_i \in k[\mathcal{L}\setminus \XX'] $. 
The association $x_{m-r-1+i,m-i}\mapsto -(1/\delta)\Gamma_i$ therefore defines a ring homomorphism
$$k[\mathcal{L}]_{\delta}/(\Delta_{\delta})=(k[X'][\mathcal{L}\setminus X'])_{\delta}/(\Delta_{\delta})\simeq k[\mathcal{L}\setminus X'])_{\delta}
$$

This entails a ring isomorphism

$$\dfrac{S_{\delta}}{\Delta_{\delta} S_{\delta}}\simeq \dfrac{k[\mathcal{L}\setminus \XX']_\delta}{(I_{m-r}(\mathcal{\widetilde{L}}))k[\mathcal{L}\setminus \XX']_\delta}.$$

Thus, ${\rm dim} \ S_{\delta}/\Delta_{\delta} S_{\delta} =\dim k[\mathcal{L}]_\delta-r-\codim (I_{m-r}(\mathcal{\widetilde{L}}))_{\delta}={\rm dim } \ S_\delta-r$

Therefore, $\codim (I_{m-r}(\mathcal{L}))$ is at least  $\codim (I_{m-r}(\mathcal{\widetilde{L}}))+r={{r+1}\choose {2}}$.

\medskip

In order to show that $I_{m-r}(\mathcal{L}) $ is the homogeneous defining ideal  of the polar variety it suffices to show that the latter has codimension at most ${{r+1}\choose {2}}$.
Since the dimension of the homogeneous coordinate ring of the polar variety coincides with the rank of the Hessian matrix of $f$, it now suffices to show
that the latter is at least $\dim R - {{r+1}\choose {2}}= m^2- {{r+1}\choose {2}}- {{r+1}\choose {2}}=m^2-r(r+1)$.  

\medskip

For this, we proceed along the same line of the proof of Theorem~\ref{cloning_generic} (ii).
Namely, set $X:=\{ x_{i,j}\,|\, i+j=r+2,r+3,\ldots,2m-r \}$ and consider the set of partial derivatives  of $f$ with respect to the variables in $X$. Let $M$ denote the Jacobian matrix of these partial derivatives with respect to the variables in $X$. Observe that $M$ is a submatrix of size  $(m^2-r(r+1))\times (m^2-r(r+1))$ of the Hessian matrix. We will show that $\det(M)\neq 0$.

Set  ${\bf v}:=\left\lbrace x_{1,m}, x_{2,m-1},\ldots,x_{m,1} \right\rbrace\subset X$, the set of variables along the main anti-diagonal of $\mathcal{DG}(r)$. 

As already pointed out, the partial derivative of $f$ with respect to any $x_{i,j}\in X $ coincides with the signed cofactor of $x_{i,j}$.  By expanding the cofactor of an entry in the set ${\bf v}$ one sees that there is a unique (nonzero) term whose support lies in ${\bf v}$ and the remaining terms have degree $\geq 2$ in the variables off ${\bf v}$. 
Similarly, the cofactor of a variable outside ${\bf v}$ has no term whose support lies in ${\bf v}$ and has exactly one (nonzero) term of degree 1 in the variables off  $\vv$. In fact, if $i+j\neq m+1$, one finds
\begin{eqnarray*}
\Delta_{j,i} &=& x_{m+1-j,m+1-i}(x_{1,m}\cdots \widehat{x_{i,m-i+1}}\cdots \widehat{x_{m-j+1,j}} \cdots x_{m,1})\\
& & \mbox{} + {\rm terms\; of\; degree\; at \;least\; 2\; off \;} \mathbf{v},
\end{eqnarray*}
where the term inside the parenthesis has support in $\mathbf{v}$.

Consider the ring endomorphism $\varphi$ of $R$ that maps any variable in {\bf v}  to itself and  any variable off {\bf v} to zero.
By the preceding observation, applying $\varphi$ to any second partial derivative of $f$ involving only the variables of $X$ will return zero or a monomial supported on the variables in ${\bf v}$. 
Let $\widetilde{M}$ denote the resulting specialized matrix of $M$.
Thus, any of its entries is either zero or a monomial supported on the variables in {\bf v}. 

We will  show that $\det(\widetilde{M}) $ is nonzero.
For this, consider the Jacobian matrix of the set of partial derivatives $\{f_v:v\in {\bf v}\}$ with respect to the variables in ${\bf v}$. Let $M_0$ denote the specialization of this Jacobian matrix by $\varphi$ considered as a corresponding submatrix of $\widetilde{M}$. Up to permutation of rows and columns of $\widetilde{M}$, we may write

$$\widetilde{M}=
\left(
\begin{array}{cc}
M_0 & N_0 \\
N_1 & M_1
\end{array}
\right),
$$
for suitable $M_1$. Now, by the way the second partial derivatives of $f$ specialize via $\varphi$ as explained above, one must have $ N_0 = N_1 = 0$. Therefore, $\det(\widetilde{M}) = \det(M_0) \det(M_1)$, so it remains to prove the nonvanishing of these two subdeterminants. Now the first block is the Hessian matrix of the form $g$ being taken as the product
of the entries in the main anti-diagonal of the matrix $\mathcal{DG}(r)$. By a similar argument used in the proof of Theorem~\ref{cloning_generic} (ii), one has that $g$ is a well-known homaloidal polynomial, hence we are done for the first matrix block. As for the second block, by construction it has exactly one nonzero entry on each row and each column. Therefore, it has a nonzero determinant.

To conclude the assertion of this item it remains to argue that the ladder determinantal ring in question is Gorenstein.
For this we use the criterion in \cite[Theorem, (b) p. 120]{ladder2}.
By the latter, we only need to see that the inner corners of the ladder depicted in Figure~\ref{r-stair-like} have indices $(a,b)$ satisfying the equality $a+b=m-1+(m-r)-1=2m-r-2$, where the ladder is a structure in an $(m-1)\times (m-1)$ matrix.

This completes the proof of this item.
The supplementary assertion on the analytic spread of $J$ is clear since the dimension of the latter equals the dimension of the $k$-subalgebra generated by the partial derivatives.
\qed

\begin{Remark}\rm
	One notes that the codimension of the polar variety in its embedding coincides with the codimension of $\mathcal{DG}(r)$ in the fully generic matrix of the same size, viewed as vector spaces of matrices over the ground field $k$.
\end{Remark}

\subsection{The ideal of the submaximal minors}

We will need a couple of lemmas, the first of which is a non-generic version of \cite[Theorem 10.16 (b)]{B-V}:

\begin{Lemma}\label{subminors_spread}
	Let $M$ be a square matrix with entries either variables over a field $k$ or zeros, such that $\det(M)\neq 0$. Let $R$ denote the polynomial ring over $k$  on the nonzero entries of $M$ and let $S\subset R$ denote the $k$-subalgebra generated by the submaximal minors.
	Then the extension $S\subset R$ is algebraic at the level of the respective fields of fractions.
\end{Lemma}

The proof is the same as the one given in \cite[Theorem 10.16 (b)]{B-V}.

The second lemma was communicated to us by Aldo Conca, as a particular case of a more general setup:

\begin{Lemma}\label{anti-diagonals}
The submaximal minors of the generic square matrix are a Gr\"obner base in the reverse lexicographic order and the initial ideal of any minor is the product of its entries along the main anti-diagonal.
\end{Lemma}

This result is the counterpart of the classical result in the case of the lexicographic order, where the initial ideals are the products of the entries along the main diagonals.
In both versions, the chosen term order should respect the rows and columns of $M$.

The content of the third lemma does not seem to have been noted before:
\begin{Lemma}\label{max-reg-seq}
Let $\mathcal{G}$ denote a generic $m\times m$ matrix and let $\mathcal{X}$ denote the set of entries none of which  belongs to the main anti-diagonal of a submaximal minor.
Then $\mathcal{X}$ is a regular sequence modulo the ideal generated by the submaximal minors in the polynomial ring of the entries of $\mathcal{G}$ over a field $k$. 
\end{Lemma}
\demo 
As for easy visualization, $\mathcal{X}$ is the set of bulleted entries below (for $m\geq 6$):
	\begin{equation}\nonumber
	\left(
	\begin{array}{ccccccccc}
	\bullet &\bullet & \bullet& \ldots &\bullet &\bullet  & \bullet  &  x_{1,m-1} & x_{1,m}\\
	\bullet & \bullet & \bullet & \ldots & \bullet & \bullet &  x_{2,m-2}&  x_{2,m-1}&  x_{2,m}\\
	\bullet & \bullet & \bullet &\ldots & \bullet & x_{3,m-3}& x_{3,m-2} & x_{3,m-1}   & \bullet\\
	\bullet  &\bullet & \bullet &\ldots & x_{4,m-4} & x_{4,m-3}& x_{4,m-2} & \bullet &\bullet  \\ [-2pt]
	\vdots &  \vdots  &\vdots & & \vdots  & \vdots & \vdots&\vdots  & \vdots\\ [-2pt]
	\bullet &x_{m-2,2} &x_{m-2,3}&\ldots & x_{m-2,m-4}&  \bullet & \bullet & \bullet & \bullet   \\[4pt]
	x_{m-1,1} & x_{m-1,2} &x_{m-1,3} &\ldots & \bullet &  \bullet & \bullet & \bullet & \bullet \\[4pt]
	x_{m,1} & x_{m,2} &\bullet & \ldots & \bullet& \bullet & \bullet & \bullet & \bullet
	\end{array}
	\right)
	\end{equation}
(A similar picture can be depicted for $m\leq 5$).

Clearly, the cardinality of $\mathcal{X}$ is $2{{m-1}\choose {2}}=(m-1)(m-2)$.	
Fix an ordering of the elements $\{a_1,\ldots,a_{(m-1)(m-2)}\}$ of $\mathcal{X}$.
By Lemma~\ref{anti-diagonals} and the assumption that every $a_i$ avoids the initial ideal of any submaximal minor, it follows that the initial ideal of the ideal $(a_1,\ldots,a_i, \mathcal{P})$ is $(a_1,\ldots,a_i, {\rm in}(\mathcal{P}))$.
Clearly, $a_{i+1}$ is not a zero divisor modulo the latter ideal, and hence, by a well known procedure, it is neither a zero divisor modulo $(a_1,\ldots,a_i, \mathcal{P})$.
\qed
 
\smallskip

In the subsequent parts we will relate the gradient ideal $J\subset R$ of the determinant of the matrix $\mathcal{DG}(r)$ in {\rm (\ref{generic-zeros})} to the ideal $I_{m-1}(\mathcal{GC})\subset R$ of its submaximal minors.
As an easy preliminary, we observe that, for any nonzero entry $x_{i,j}$ of $\mathcal{DG}(r)$, since the nonzero entries of the matrix are independent variables, it follows easily from the Laplace expansion along the $i$th row that the $x_{i,j}$-derivative $f_{i,j}$ of $f$ coincides with the (signed) cofactor of  $x_{i,j}$. In particular, one has $J\subset I_{m-1}(\mathcal{GC})$ throughout the entire subsequent discussion and understanding the conductor  $J: I_{m-1}(\mathcal{GC})$ will be crucial.

\begin{Proposition}\label{sizing_the_conductor}
Let $\mathcal{DG}(r)$ as in {\rm (\ref{generic-zeros})} denote our basic degenerate matrix, with $1\leq r\leq m-2$.
For every $0\leq j\leq m$, consider the  submatrices $M_j$ and $N_j$ of $\mathcal{DG}(r)$  consisting of its last $j$  columns and the  its last  $j$ rows, respectively.
Write $I:=I_{m-1}(\mathcal{GC})\subset R$ for the ideal of $(m-1)$-minors of $\mathcal{DG}(r)$ and $J$ for the gradient ideal of $f:=\det(\mathcal{DG}(r))$.
Then 
$I_j(N_{j})\cdot I_{r-j}(M_{r-j})\subset J:I$ for every $0\leq j\leq r.$
\end{Proposition}
\demo
For a fixed  $1\leq j\leq r,$ we write the matrices $\mathcal{DG}(r)$ and its adjoint adj$(\mathcal{DG}(r))$  in  the  following block form:

\begin{equation}
\mathcal{DG}(r)=\left(\begin{array}{cccc}
\widetilde{N}_j\\ [3pt]
N_j
\end{array}\right),
\quad
{\rm adj}(\mathcal{DG}(r))=\left(\begin{array}{c|ccc}
\Theta_{1,j}&\Theta_{2,j}\\
\hline
\Theta_{3,j}&\Theta_{4,j}
\end{array}\right);
\end{equation}
where $\Theta_{1,j},\Theta_{2,j},\Theta_{3,j},\Theta_{4,j}$
stand for submatrices of sizes 
$(j+m-r)\times(m-j),  (j+m-r)\times j,  (r-j)\times(m-j)$ and $(r-j)\times j$, respectively.
Thus, we have

\begin{equation}\label{13}
{\rm adj}(\mathcal{DG}(r))\cdot\mathcal{DG}(r)=\left(\begin{array}{cccc}
\Theta_{1,j}\widetilde{N}_j+\Theta_{2,j}N_j\\ [3pt]
\Theta_{3,j}\widetilde{N}_j+\Theta_{4,j}N_j
\end{array}\right)=f\cdot \mathbb{I}_m.
\end{equation}
with $\mathbb{I}_m$ denoting the identity matrix of order $m$.
Since  $f$ belongs to $J,$ then $I_1(\Theta_{1,j}\widetilde{N}_j+\Theta_{2,j}N_j)\subset J$. On the other hand, the entries of  $\Theta_{1,j}$ are cofactors of the entries on the upper left corner of $\mathcal{DG}(r)$, hence belong to  $J$ as well. 
Therefore $I_1(\Theta_{2,j}N_j)\subset J$ as well.
From this by an easy argument it follows that
\begin{equation}
\label{inclusao1}
I_1(\Theta_{2,j})I_j(N_j)\subset J
\end{equation}
and, for even more reason,
\begin{equation}\label{inclusao2}
\xymatrix{*+[F]{I_1(\Theta_{2,j})I_j(N_j)\cdot I_{r-j}(M_{r-j})\subset J}}
\end{equation}

Similarly, writing
$$\mathcal{DG}(r)=\left(\begin{array}{c|ccc}
\widetilde{M}_{r-j}&
M_{r-j}
\end{array}\right)
$$ 
we have:
{\small
	\begin{equation*}
	\mathcal{DG}(r)\cdot {\rm adj}(\mathcal{DG}(r))=\left(\begin{array}{c|ccc}
	\widetilde{M}_{r-j}\Theta_{1,j}+M_{r-j}\Theta_{3,j}&\widetilde{M}_{r-j}\Theta_{2,j}+M_{r-j}\Theta_{4,j}
	\end{array}\right)=f\cdot {\rm Id}_m.
	\end{equation*}}

An entirely analogous reasoning leads to the inclusion
$I_1(\Theta_{3,j})I_{r-j}(M_{r-j})\subset J$, and for even more reason
\begin{equation}\label{inclusao4}
\xymatrix{*+[F]{I_1(\Theta_{3,j})I_j(N_j)\cdot I_{r-j}(M_{r-j})\subset J}.}
\end{equation} 

Arguing now with the second block  $\widetilde{M}_{r-j}\Theta_{2,j}+M_{r-j}\Theta_{4,j}$, again $I_1(\widetilde{M}_{r-j}\Theta_{2,j}+M_{r-j}\Theta_{4,j})\subset J$, and hence  for each $\delta\in I_j(N_j),$  also  $I_1(\delta\widetilde{M}_{r-j}\Theta_{2,j}+\delta M_{r-j}\Theta_{4,j})\subset J.$ But, by \eqref{inclusao1}, the entries of $\delta\widetilde{M}_{r-j}\Theta_{2,j}$ belong to  $J.$ Thus, the entries de $\delta M_{r-j}\Theta_{4,j}$ belong to $J$ and consequently
\begin{equation}\label{inclusao5}
\xymatrix{*+[F]{I_1(\Theta_{4,j})I_j(N_j)\cdot I_{r-j}(M_{r-j})\subset J}}.
\end{equation}

It follows from \eqref{inclusao2}, \eqref{inclusao4} and \eqref{inclusao5}  that
\begin{equation}
\left(I_1(\Theta_{2,j}),\,I_1(\Theta_{3,j}),\,I_1(\Theta_{4,j})\right)\,I_j(N_j)\cdot I_{r-j}(M_{r-j}))\subset J.
\end{equation}
Since also $I_1(\Theta_{1,j})\subset J,$  we have 
\begin{equation}\label{inclusao6}
(I_1(\Theta_{1,j}),\,I_1(\Theta_{2,j}),\,I_1(\Theta_{3,j}),\,I_1(\Theta_{4,j}))\,I_j(N_j)\cdot I_{r-j}(M_{r-j}))\subset J.
\end{equation}
From this equality it obtains 
\begin{equation}\label{inclusao7}
I_j(N_j)\cdot I_{r-j}(M_{r-j})I\subset J 
\end{equation}
because
$$I=I_{m-1}(\mathcal{DG}(r))=I_1({\rm adj}(\mathcal{DG}(r)))=(I_1(\Theta_{1,j}),\,I_1(\Theta_{2,j}),\,I_1(\Theta_{3,j}),\,I_1(\Theta_{4,j})).$$
This establishes the assertion above -- we note that it contains as a special case (with $j=0$ and $j=r$) the separate inclusions $I_r(M_r),\,I_r(N_r)\subset J:I$.
\qed

\begin{Theorem}\label{primality_generic_zeros} 
	Consider the matrix $\mathcal{DG}(r)$ as in {\rm (\ref{generic-zeros})}, with $1\leq r\leq m-2$.
	Let $I:=I_{m-1}(\mathcal{DG}(r))\subset R$ denote its ideal of $(m-1)$-minors and $J$ the gradient ideal of $f:=\det(\mathcal{DG}(r))$.
	Then
\begin{enumerate}
\item[{\rm (i)}] $I$ is a Gorenstein ideal of codimension $4$ and maximal analytic spread.
\item[{\rm (ii)}] The $(m-1)$-minors of $\mathcal{DG}(r)$ define a birational map $\pp^{m^2-{{r+1}\choose {2}}-1}\dasharrow \pp^{m^2-1}$ onto a cone over the polar variety of $f$ with vertex cut by ${{r+1}\choose {2}}$ coordinate hyperplanes.
\item[{\rm (iii)}] The conductor $J:I$ has codimension  $2(m-r)\geq 4$; in particular, $J$ has codimension $4$. 
\item[{\rm (iv)}]  If $r\leq m-3$ then $I$ is contained in the unmixed part of $J;$  in particular,  if $R/J$ is Cohen--Macaulay then $r=m-2$. 
\item[{\rm (v)}]  Let $r=m-2$. Then the set of minimal primes of $J$ is exactly the set of associated primes of $I$ and  of the ideals $(I_{j}(N_j),\,I_{m-1-j}(M_{m-1-j}))$, for $1\leq j\leq m-2$.
\item[{\rm (vi)}]  If, moreover, ${{r+1}\choose {2}}\leq m-3$, then $I$ is a prime ideal$;$ in particular, in this case it coincides with the unmixed part of $J$. 
	\end{enumerate}
\end{Theorem}
\demo
(i) The analytic spread follows from Lemma~\ref{subminors_spread}.

The remaining assertions of the item follow from Lemma~\ref{max-reg-seq}, which shows that $I$ is a specialization of the ideal of generic submaximal minors, provided we argue that the set 
{\large 
	$$\{x_{m,m-r+1}; x_{m,m-r+2}, x_{m-1,m-r+2}; \ldots; x_{m,m}, x_{m-1,m},\ldots, x_{m-r+1,m}\}
	$$}
of variables on the voided entry places of the generic $m\times m$ matrix
{\small
	\begin{equation}\nonumber
	\left(
	\begin{array}{cccccccc}
	x_{1,1}&\ldots & x_{1,m-r}& x_{1,m-r+1} &  x_{1,m-r+2}&\ldots & x_{1,m-1} & x_{1,m}\\
	\vdots & \ldots & \vdots & \vdots & \vdots & \ldots & \vdots &\vdots \\
	x_{m-r,1}& \ldots & x_{m-r,m-r}& x_{m-r,m-r+1} & x_{m-r,m-r+2}&\ldots & x_{m-r,m-1} & x_{m-r,m}\\
	x_{m-r+1,1}&\ldots & x_{m-r+1,m-r}& x_{m-r+1,m-r+1} & x_{m-r+1,m-r+2} &\ldots & x_{m-r+1,m-1} & \\
	x_{m-r+2,1}&\ldots & x_{m-r+2,m-r}& x_{m-r+2,m-r+1} & x_{m-r+2,m-r+2} &\ldots &  & \\
	\vdots & \ldots & \vdots & \vdots & \vdots&\iddots &\vdots & \vdots\\
	x_{m-1,1} &\ldots & x_{m-1,m-r}&  x_{m-1,m-r+1}& &\ldots &  & \\[4pt]
	x_{m,1} &\ldots & x_{m,m-r}& & &\ldots &  & \\ [3pt]
	\end{array}
	\right)
	\end{equation}}
is a subset of $\mathcal{X}$ as in the lemma.
But this is immediate  because of the assumption $r\leq m-2$.

\smallskip

(ii) By Lemma~\ref{max-reg-seq}, the ideal is a specialization of the ideal of submaximal minors in the generic case; in particular, it is linearly presented.
On the other hand, its analytic spread is maximal by Lemma~\ref{subminors_spread}.
Therefore, by Theorem~\ref{basic_criterion} the minors define a birational map onto the image. It remains to argue that the image is a cone over the polar variety, the latter as described in Theorem~\ref{polar-zeros} (c).

To see this note the homogeneous inclusion $T:=k[J_{m-1}]\subset T':=k[I_{m-1}]$ of $k$-algebras which are domains, where $I$ is minimally generated by the generators of $J$ and by ${{r+1}\choose {2}}$ additional generators, say, $f_1,\ldots, f_s$, where $s={{r+1}\choose {2}}$, that is, $T'=T[f_1,\ldots, f_s]$.
On the other hand, one has $\dim T=m^2-r(r+1)$ and  $\dim T'=m^2-{{r+1}\choose {2}}$.
Therefore, ${\rm tr. deg}_{k(T)} k(T)(f_1,\ldots, f_s)=\dim T'-\dim T= {{r+1}\choose {2}}=s$, where $k(T)$ denotes the field of fractions of $T$.
This means that $f_1,\ldots, f_s$ are algebraically independent over $k(T)$ and, a fortiori, over $T$.
This shows that $T'$ is a polynomial ring over $T$ in ${{r+1}\choose {2}}$ indeterminates. Geometrically, the image of the map defined by the $(m-1)$-minors is a cone over the polar image with vertex cut by ${{r+1}\choose {2}}$ independent hyperplanes.

\smallskip

(iii) 
The fact that the stated value $2(m-r)$ is an upper bound follows from the following fact: first, because the only cofactors which are not partial derivatives (up to sign) are those corresponding to the zero entries, then $J$ is contained in the ideal $Q$ generated by the variables of the last row and the last column of the matrix. 
Since $Q$ is prime and $I\not\subset Q$ then clearly $J:I\subset Q$.

To see that $2(m-r)$ is a lower bound as well, we will use Proposition~\ref{sizing_the_conductor}.

For that we need some intermediate results.

\smallskip

{\sc Claim 1.} For every $1\leq j\leq r$, both $I_j(M_j)$ and $I_j(N_j)$ have codimension $m-r$.

By the clear symmetry, it suffices to consider $I_j(M_j)$.
Note that $M_j$ has $r-j+1$ null rows, so its ideal of $j$-minors coincides with the ideal of $j$-minors of its $(m-(r-j+1))\times j$ submatrix  $M_j'$ with no null rows.
Clearly, this ideal of (maximal) minors has codimension at most $(m-(r-j+1))-j+1= m-r$.  Now,  the matrix  $M_j'$  specializes to the well-known diagonal specialization using only $m-r$ variables -- by definition, the latter is the specialization of a suitable Hankel  matrix via the ring homomorphism  mapping to zero the variables  of the upper left and lower right corner except the last variables of first column and  the first variable of the last column. This ensures that $I_j(M_j')$ has codimension at least $m-r$.

\smallskip

{\sc Claim 2.} For every $1\leq j\leq r-1,$ the respective sets of nonzero entries of $N_{j}$ and $M_{r-j+1}$ are disjoint. In particular, the codimension of $(I_{j}(N_j),\,I_{r-j+1}(M_{r-j+1}))$ is $2(m-r).$

The disjointness assertion is clear by inspection and the codimension follows from the previous claim.

\smallskip

To proceed, we envisage the following chains of inclusions

\begin{equation}
I_1(N_1)\supset I_2(N_2)\supset\ldots\supset I_{r-1}(N_{r-1})\supset I_{r}(N_r)
\end{equation}
and
\begin{equation}
I_1(M_1)\supset I_2(M_2)\supset\ldots\supset I_{r-1}(M_{r-1})\supset I_{r}(M_r).
\end{equation}

Let $P$ denote a prime ideal containing the conductor $J:I.$
By Proposition~\ref{sizing_the_conductor} one has the inclusion  $I_1(N)\cdot I_{r-1}(M_{r-1})\subset P.$ Thus,

\begin{enumerate}
	\item[($A_1$)] either $I_{1}(N_1)\subset P$, or else
	\item[($B_1$)] $I_1(N_1)\not\subset J:I$ but $I_{r-1}(M_{r-1})\subset P$.
\end{enumerate}

If $(A_1)$ is the case, then  $(I_1(N_1),\,I_r(M_r))\subset P$,  because $I_r(M_r)\subset J:I$ again by Proposition~\ref{sizing_the_conductor} (with $j=0$).
By Claim 2 above, we then see that the codimension of $J:I$ is at least $2(m-r)$.

If $(B_1)$ takes place then we consider the inclusion  $I_2(N_2)\cdot I_{r-2}(M_{r-2})\subset J:I\subset P$ by Proposition~\ref{sizing_the_conductor}. 
The latter in turn gives rise to two possibilities according to which
\begin{enumerate}
	\item[($A_2$)] either $I_{2}(N_2)\subset P,$ or else
	\item[($B_2$)] $I_2(N_2)\not\subset P$ but $I_{r-2}(M_{r-2})\subset P.$
\end{enumerate}

Again, if  $(A_2)$ is the case then $(I_2(N_2),\,I_{r-1}(M_{r-1}))\subset P$ since $I_{r-1}(M_{r-1})\subset P$ by hypothesis. Once more, by Claim 2, the codimension of $P$ is at least $2(m-r).$

If intead $(B_2)$ occurs then we step up to the inclusion $I_3(N_3)\cdot I_{r-3}(M_{r-3})\subset P$ and repeat the argument. 
Proceeding in this way, we may eventually find an index $1\leq j\leq r-1$ such that the first alternative $(A_j)$ holds, in which case we are through always by Claim 2.  Otherwise, we must be facing the situation where  $I_{j}(N_{j})\not\subset P $ for every $1\leq j\leq r-1.$ In particular, $I_{r-1}(N_{r-1})\not\subset P$ and $I_1(M_1)\subset P.$ Thus, 
$(I_r(N_r),\,I_1(M_1))\subset P,$ and once more by Claim 2, $P$ has codimension at least $2(m-r).$ This concludes the proof of the codimension of $J:I$.

The assertion  that $J$ has codimension $4$ is then ensured as $J\subset I$ and $I$ has codimension $4$ by item (i). 	
\qed

\smallskip

(iv) By (iii), if $r\leq m-3$ then $J:I$ has codimension at least $2(m-r)\geq 6$. This implies that $I\subset J^{\rm un}$ and the two coincide up to radical.

The first assertion on the Cohen--Macaulayness of $R/J$ is clear since then $J$ is already unmixed, hence $J=I$ which is impossible since $r\geq 1$.

\smallskip

(v) First, one has an inclusion $J\subset (I_{j}(N_j),\,I_{m-1-j}(M_{m-1-j}))$.
Indeed, recall once more that the partial derivatives are (signed) cofactors.
Expanding each cofactor by Laplace along $j$-minors, one sees that it is expressed as products of generators of $I_{j}(N_j)$ and of $I_{m-1-j}(M_{m-1-j})$.
Next, the ideal $ (I_{j}(N_j),\,I_{m-1-j}(M_{m-1-j}))$ is perfect of codimension $4$ -- the codimension is clear by Claim 2 in the proof of (iii), while perfectness comes from the fact that the tensor product over $k$ of Cohen--Macaulay $k$-algebras of finite type is Cohen--Macaulay.

At the other end, by (i) the ideal $I$ is certainly perfect.
Therefore, any associated prime of either $(I_{j}(N_j),\,I_{m-1-j}(M_{m-1-j}))$  or $I$ is a minimal prime of $J$.

Conversely, let $P$ denote a minimal prime of $J$ not containing $I$. Then $J:I\subset P$, hence the same argument in the proof of (iii) says that $P$ contains some ideal of the form  $(I_{j}(N_j),\,I_{m-1-j}(M_{m-1-j}))$.

\smallskip

(vi)  We will apply Proposition~\ref{Pisprime} in the case where $M'=\mathcal{G}$ is an $m\times m$ generic matrix  and $M=\mathcal{DG}(r)$ is the degenerated generic matrix as in the statement.
In addition, we take $k=m-2$, so $k+1=m-1$ is the size of the submaximal minors.
Observe that the vector space spanned by the entries of $M$ has codimension ${{r+1}\choose {2}}$ in the  vector space spanned by the entries of $M'$. Since  the $m\times m$ generic matrix is  $2=m-(m-2)$-generic (it is $m$-generic as explained in \cite[Examples, p. 548]{Eisenbud2}), the theorem ensures that if  ${{r+1}\choose {2}}\leq k-1=m-3$ then  $I_{m-1}(\mathcal{DG})$ is  prime.

Since in particular $r\leq m-3$ then item (d) says that $I\subset J^{\rm un}$.
But $I$ is prime, hence $I= J^{\rm un}$.
\qed

\begin{Remark}\rm
	The statement of item (i) in Theorem~\ref{primality_generic_zeros} depends not only the number of the entries forming a regular sequence on $\mathcal{P}$ but also their mutual position. Thus, for example, if more than $r$ of the entries belong to one same column or row it may happen that $I$ has codimension strictly less than $4$.
\end{Remark}

We end by filing a few natural questions/conjectures.
The notation is the same as in the last theorem.

\begin{Question}
	Does  $I=J^{\rm un}$ hold for $r\leq m-3$? {\rm (}By {\rm (iv)} above it would suffice to prove that $I$ is a radical ideal.{\rm )}
\end{Question}

\begin{Question}
	Is ${{r+1}\choose {2}}\leq m-3$ the exact obstrution for the primality of the ideal $I$ of submaximal minors?
\end{Question}

\begin{Conjecture}
	If ${{r+1}\choose {2}}\leq m-3$ then $J$ has no embedded primes.
\end{Conjecture}

\begin{Conjecture} Let $r=m-2$. Then for any $1\leq j\leq m-2$ one has the following primary decomposition:
	$$I_j(M_j)=(x_{j+1,m-j+1},\delta_{j})\cap(x_{j,m-j+2},\delta_{i-1})\cap\ldots\cap(x_{3,m-1},\delta_{2})\cap(x_{2,m},x_{1,m}),$$
	where $\delta_{t}$ denotes the determinant of the $t\times t$ upper submatrix of $M_t$.
	A similar result holds for $I_j(N_j)$ upon reverting the indices of the entries and replacing $\delta_{t}$ by the determinant $\gamma_t$ of the $t\times t$ leftmost submatrix of $N_t$.
	In particular, $I_j(M_j)$ and $I_j(N_j)$ are radical ideals.
\end{Conjecture}

\begin{Conjecture} In the notation of the previous conjecture, one has
\begin{eqnarray}\nonumber
J:I&=&\left(\bigcap_{t=1}^{r}(x_{m,1},x_{m,2},x_{t+1,m-t+1},\delta_t)\right)
\cap\left(\bigcap_{t=1}^{r-1}(x_{m-1,3},\gamma_{2},x_{t+1,m-t+1},\delta_t)\right) \cap\\ \nonumber
 &\ldots &\cap\left(\bigcap_{t=1}^{1}(x_{m,1},x_{m,2},x_{t+1,m-t+1},\delta_t)\right)
 \cap\left(x_{3,m-1},\gamma_{r},x_{2,m},\delta_1\right).
\end{eqnarray}
In particular, $J:I$ is a radical ideal.
\end{Conjecture}

\begin{Conjecture}
If $r=m-2$ then both $R/J$ and $R/J:I$ are Cohen--Macaulay reduced rings and, moreover, one has $J=I\cap (J:I)$.
\end{Conjecture}

\subsection{The dual variety}

We keep the previous notation with $\mathcal{DG}(r)$ denoting the $m\times m$ matrix in (\ref{generic-zeros}).

In this part we describe the structure of the dual variety $V(f)^*$ of $V(f)$ for $f=\det \mathcal{DG}(r)$. 
The result in particular answers affirmatively a question posed by F. Russo as to whether the codimension of the dual variety of a homogeneous polynomial in its polar variety can be arbitrarily large when its Hessian determinant vanishes.
In addition it shows that this can happen in the case of structured varieties.

As a preliminary, we file the following lemma which may have independent interest.

\begin{Lemma}\label{specializing_max_minors} Let $\mathcal{G}_{a,b}$ denote the $a\times b$ generic matrix, with $a\geq b$. Letting $0\leq r\leq b-2$, consider the corresponding degeneration matrix
$\Psi:=\mathcal{DG}_{a,b}(r)$:
	
	$$\left(\begin{array}{cccccccc}
	x_{1,1}&x_{1,2}&\ldots&x_{1,b-r}&x_{1,b-r+1}&\ldots&x_{1,b-1}&x_{1,b}\\
	x_{2,1}&x_{2,2}&\ldots&x_{2,b-r}&x_{2,b-r+1}&\ldots&x_{2,b-1}&x_{2,b}\\
	\vdots&\vdots&&\vdots&\vdots&&\vdots&\vdots\\
	x_{a-r,1}&x_{a-r,2}&\ldots&x_{a-r,b-r}&x_{a-r,b-r+1}&\ldots&x_{a-r,b-1}&x_{a-r,b}\\
	x_{a-r+1,1}&x_{a-r+1,2}&\ldots&x_{a-r+1,b-r}&x_{a-r+1,b-r+1}&\ldots&x_{a-r+1,b-1}&0\\
	\vdots&\vdots&&\vdots&\vdots&\iddots&\vdots&\vdots\\
	x_{a-1,1}&x_{a-1,2}&\ldots&x_{a-1,b-r}&x_{a-1,b-r+1}&\ldots&0&0\\
	x_{a,1}&x_{a,2}&\ldots&x_{a,b-r}&0&\ldots&0&0
	\end{array}\right).$$
	Then the ideal of maximal minors of $\Psi$ has the expected codimension $a-b+1.$
\end{Lemma}
\demo The argument is pretty much the same as in the proof of Lemma~\ref{max-reg-seq}, except that the right lower corner of $\mathcal{G}_{a,b}$ whose entries will form a regular sequence is now as depicted in blue below 
$$
X=\left(\begin{array}{cccccccc}
x_{1,1}&x_{1,2}&x_{1,3}&\ldots&x_{1,b-1}&x_{1,b}\\
x_{2,1}&x_{2,2}&x_{2,3}&\ldots&x_{2,b-1}&x_{2,b}\\
\vdots&\vdots&\vdots&&\vdots&\vdots\\
x_{a-b+2,1}&x_{a-b+2,2}&x_{a-b+2,3}&\ldots&x_{a-b+2,b-1}&\textcolor{blue}{x_{a-b+2,b}}\\
x_{a-b+3,1}&x_{a-b+3,2}&x_{a-b+3,3}&\ldots&\textcolor{blue}{x_{a-b+3,b-1}}&\textcolor{blue}{x_{a-b+3,b}}\\
\vdots&\vdots&\vdots&\iddots&\vdots&\vdots\\
x_{a-1,1}&x_{a-1,2}&\textcolor{blue}{x_{a-1,3}}&\ldots&\textcolor{blue}{x_{a-1,b-1}}&\textcolor{blue}{x_{a-1,b}}\\
x_{a,1}&\textcolor{blue}{x_{a,2}}&\textcolor{blue}{x_{a,3}}&\ldots&\textcolor{blue}{x_{a,b-1}}&\textcolor{blue}{x_{a,b}}\\
\end{array}\right).$$
(Note that the ``worst'' case is the $b$-minor with of the last $b$ rows of $\mathcal{G}_{a,b}$, hence the top right blue entry above.)

To finish the proof of the lemma, one argues that the blue entries form a regular sequence on the initial ideal of $I_b(\mathcal{G}_{a,b})$ since the latter is generated by the products of the entries along the anti-diagonals of the maximal minors (in any monomial order).
\qed

\medskip

Next is the main result of this part. We stress that, in contrast to the dual variety in the case of the cloning degeneration, here the dual variety will in fact be a ladder determinantal variety. 

Back to the notation of (\ref{generic-zeros}), one has:

\begin{Theorem}\label{dual_zeros} Let $0\leq r\leq m-2$ and $f:=\det \mathcal{DG}(r)$. Then:
	\begin{enumerate}	
		\item[{\rm (a)}] The dual variety $V(f)^*$ of $V(f)$ is a ladder determinantal variety of codimension $(m-1)^2 - {r+1 \choose 2}$ defined by $2$-minors; in particular it is  arithmetically Cohen--Macaulay and its codimension  in the polar variety of $V(f)$ is $(m-1)^2-r(r+1)$.
		\item[{\rm(b)}] $V(f)^*$  is arithmetically Gorenstein if and only if $r=m-2$.
	\end{enumerate}
\end{Theorem}
\demo
(a)
The proof will proceed in parallel to the proof of Theorem~\ref{dim_dual}, but there will be some major changes.

We first show that once more $\dim V(f)^*=2m-2$ and, for that, check the side-inequalities separately.

{\bf 1.} $\dim V(f)^*\geq  2m-2$. 

At the outset we draw as before upon the equality of Segre (\cite{Segre}):
$$\dim V(f)^*=\rk H(f) \pmod{f}-2,$$ 
where $H(f)$ denotes the Hessian matrix of $f$.
It will then suffice to show that $H(f)$ has a submatrix of rank at least $2m$ modulo $f$.

For this purpose, consider the submatrix $\phi$ of $\mathcal{DG}(r)$ obtained by omitting the first row:
{\small
	\begin{equation*}\label{generic-zeros}
	\left(
	\begin{array}{cccccccc}
	x_{2,1}&\ldots & x_{2,m-r}& x_{2,m-r+1} &  x_{2,m-r+2}&\ldots & x_{2,m-1} & x_{2,m}\\
	\vdots & \ldots & \vdots & \vdots & \vdots & \ldots & \vdots &\vdots \\
	x_{m-r,1}& \ldots & x_{m-r,m-r}& x_{m-r,m-r+1} & x_{m-r,m-r+2}&\ldots & x_{m-r,m-1} & x_{m-r,m}\\
	x_{m-r+1,1}&\ldots & x_{m-r+1,m-r}& x_{m-r+1,m-r+1} & x_{m-r+1,m-r+2} &\ldots & x_{m-r+1,m-1} & 0\\
	x_{m-r+2,1}&\ldots & x_{m-r+2,m-r}& x_{m-r+2,m-r+1} & x_{m-r+2,m-r+2} &\ldots & 0 & 0\\
	\vdots & \ldots & \vdots & \vdots & \vdots&\iddots &\vdots & \vdots\\
	x_{m-1,1} &\ldots & x_{m-1,m-r}&  x_{m-1,m-r+1}& 0&\ldots & 0 & 0\\[4pt]
	x_{m,1} &\ldots & x_{m,m-r}& 0& 0&\ldots & 0 & 0\\ [3pt]
	\end{array}
	\right)
	\end{equation*}
}

 {\sc Claim:} The ideal $I_t(\phi)$ has codimension at leat $m-t+1$ for every $1\leq t\leq m-1$.
 
To see this, for every $1\leq t\leq m-1,$ consider the $m-1\times t$ submatrix $\Psi_t$ of $\phi$ of the first $t$ columns. 

By Lemma~\ref{specializing_max_minors}, the codimension of $I_t(\Psi_t)$ is at least $m-t$ -- note that the lemma is applied with $a=m-1, b=t\leq m-1$ and $r$ updated to $r':=r-(m-t)$, so indeed $b-r'=t-r+(m-t)=m-r\geq 2$ as required.

It remains to find some $t$-minor of $\phi$ which is a nonzerodivisor on $I_t(\Psi_t)$.
One choice is the $t$-minor $D$ of the columns  $2,3\ldots, t-1, m$ and rows $1, m-t+1,m-t+2,\ldots,m-1$.
Indeed, a direct verification shows that the variables in the support of the initial ${\rm in}(D)$ of $D$ in the reverse lexicographic order do not appear in the supports of the initial terms of the maximal minors of $\Psi_t$. Therefore, ${\rm in}(D)$ is a nonzerodivisor on the initial ideal of  $I_t(\Psi_t)$. A standard argument then shows that  $I_t(\Psi_t)$ has codimension at least $m-t+1$.

This completes the proof of the claim.

The statement means that the ideal of maximal minors of $\phi$ satisfies the so-called property $(F_1)$. Since it is a perfect ideal of codimension $2$ with $\phi$ as its defining Hilbert--Burch syzygy matrix, it follows as in the proof of Theorem~\ref{dim_dual} that it is an ideal of linear type.
(\cite{Trento}).
In particular the $m$ maximal minors of $\phi$ are algebraically independent over $k$, hence their Jacobian matrix with respect to the entries of $\phi$ has rank $m$.

Form this point the argument proceeds exactly as in the  proof of item {\bf 1.} of Theorem~\ref{dim_dual}.

\smallskip

{\bf 2.} $\dim V(f)^*\leq  2m-2$.

Let $P\subset k[\yy]:=k[y_{i,j}\,|\, 2\leq i+j\leq 2m-r]$ denote the homogeneous defining ideal of the dual variety $V(f)^*$ in its natural embedding, i.e., 
\begin{equation}
k[\partial f/\partial x_{i,j}\,|\,2\leq i+j\leq 2m-r]/(f)\simeq k[\yy]/P.
\end{equation}
The isomorphism is an isomorphism of graded $k$-algebras induced by the assignment $y_{i,j}\mapsto \partial f/\partial x_{i,j}.$ 

\smallskip

{\sc Claim 2.} The homogeneous defining ideal $P$ of $V(f)^*$ contains the ideal generated by the $2\times 2$ minors of the following ladder matrix:


\begin{figure}[h]
	\centering
	\includegraphics[width=15.5cm, height=5.5cm]{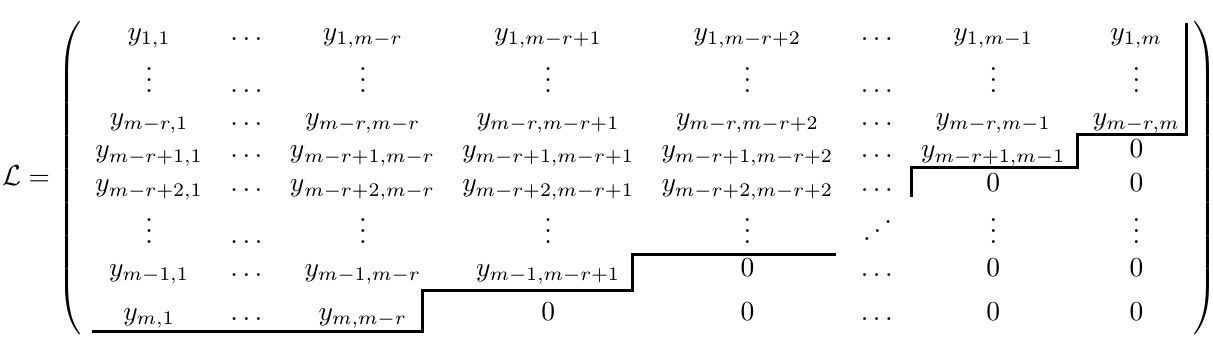}\\
	\caption{Ladder ideal of the dual variety}\label{fig3}
\end{figure}

To see  this, we first recall that by Proposition~\ref{GolMar}, $\partial f/\partial x_{i,j}$ coincides with the cofactor of $x_{i,j}$ in $\mathcal{DG}(r)$. Now, consider the following relation afforded by the cofactor identity:
\begin{equation}\label{relation_mod_f}
{\rm adj}(\mathcal{DG}(r))\cdot\mathcal{DG}(r)\equiv 0\,({\rm mod\,}f).
\end{equation}
Further, for each pair of integers $i,j$ such that $1\leq i<j\leq m $ let $F_{ij}$ denote the $2\times m$ submatrix of ${\rm adj}(\mathcal{DG}(r))$ consisting of the $i$th and $j$th rows.
In addition, let $C$ stand for the $m\times (m-1)$ submatrix of $\mathcal{DG}(r)$  consisting of its $m-1$ leftmost columns.
Then (\ref{relation_mod_f}) gives the relations
$$F_{ij}C\equiv0\,({\rm mod\,}f),$$
for all $1\leq i<j\leq m $.
From this, since the rank of $C$ modulo $(f)$ is obviously still $m-1$, the rank of every $F_{i,j}$ is necessarily $1$.
This shows that every $2\times 2$ minor of ${\rm adj}(\mathcal{GD}(r))$ vanishes modulo $(f)$.
Therefore,  each such minor involving only cofactors that are partial derivatives gives a $2\times 2$ minor of $\mathcal{L}$ vanishing  on the partial derivatives. Clearly, by construction, we obtain this way all the $2\times 2$ minors of $\mathcal{L}$.
This proves the claim.

\smallskip

Now, since $I_2(\mathcal{L})$ is a ladder determinantal ideal on a suitable  generic matrix it is a Cohen-Macaulay prime ideal (see \cite{Nar} for primeness and \cite{HeTr} for Cohen--Macaulayness). 
Moreover, its codimension is $m^2-{r+1 \choose 2}-(2m-1)=(m-1)^2 - {r+1 \choose 2}$ as follows from an application to this case of the general principle in terms of maximal chains as described in \cite[Theorem 4.6 and Corollary 4.7]{HeTr}.
\qed

\medskip

(b) By the previous item, the homogeneous defining ideal of the dual variety is generated by the $2\times 2$ minors of the ladder matrix in Figure~\ref{fig3}. Observe that  the smallest square matrix containing all the entries of the latter is the $m\times m$ matrix $\mathcal{DG}(r)$. By \cite[Theorem, (b) p. 120]{ladder2},  the ladder ideal is  Gorenstein if only if the inner corners of the ladder  have indices $(i,j)$ satisfying the equality $i+j=m+1$.
In the present case, the inner corners have indices satisfying the equation 
$$i+j=m-r+ (m-1)=m-r+1+(m-2)=\cdots = m-1 +(m-r)=2m-r-1.$$
Clearly this common value equals $m+1$ if and only if $r=m-2$.
\qed

\begin{Remark}\rm
	It may of some interest to note that degenerating the matrix $\mathcal{DG}(m-2)$ all the way to a Hankel matrix,  recovers the so-called sub-Hankel matrix thoroughly studied in \cite{CRS} from the homaloidal point of view and in \cite{MAron}, from the ideal theoretic side. The situation is ever more intriguing since in the sub-Hankel case the determinant is homaloidal.
\end{Remark}


\end{document}